\documentclass[12pt,a4paper]{amsart}
\usepackage{amssymb,amsxtra}
\usepackage[english,french]{babel}

\newtheorem{thm}{Theorem}[section]
\newtheorem{lem}[thm]{Lemma}
\newtheorem{prop}[thm]{Proposition}
\newtheorem{cor}[thm]{Corollary}

\theoremstyle{definition}
\newtheorem{defn}[thm]{Definition}

\theoremstyle{remark}
\newtheorem{remark}[thm]{Remark}

\theoremstyle{plain}

\theoremstyle{remark}

\newtheorem*{example}{Example}
\numberwithin{equation}{section}

\begin{document}

\title{ Exponential sums and Rank of Persymmetric Matrices over $\mathbb{F}_{2} $}
\author{Jorgen~Cherly}
\address{D\'epartement de Math\'ematiques, Universit\'e de
    Brest, 29238 Brest cedex~3, France}
\email{Jorgen.Cherly@univ-brest.fr}
\thanks{The author is indepted to Thierry Levasseur UBO Brest for helpful advice regarding LaTeX}

\maketitle 
\begin{abstract}
 Soit  $\mathbb{K} $ le corps des s\'eries  de Laurent formelles  $ \mathbb{F}_{2}((T^{-1})). $
Nous calculons  en particulier des sommes  exponentielles dans  $\mathbb{K} $  de la forme $ 
 \sum_{degY\leq  k-1}\sum_{degZ\leq s-1}E(tYZ) $ o\`{u}
 t  est dans la boule unit\'{e} de  $\mathbb{K},$  en d\'{e}montrant qu'elles  d\'{e}pendent  seulement du rang de matrices 
 persym\'{e}triques  avec des entr\'{e}es dans $\mathbb{F}_{2} $ qui leur sont associ\'{e}es.
  ( Une matrice \;$ [\alpha _{i,j}]  $ est   persym\'{e}trique  si  $ \alpha _{i,j} = \alpha _{r,s} $ \; pour  \; i+j = r+s ). 
 En outre nous \'{e}tablissons des propri\'{e}tes de rang d'une partition de matrices persym\'{e}triques.
 Nous utilisons ces r\'{e}sultats pour calculer le nombre $\Gamma _{i}$ de matrices  persym\'{e}triques sur $\mathbb{F}_{2}$
 de rang i.
 Nous retrouvons en particulier une formule g\'{e}n\'{e}rale donn\'{e}e par D.E.Daykin.
 Notre dŽmonstration est, comme indiqu\'{e}, tr\`{e}s  diff\'{e}rente, puisqu'elle se fonde sur les propri\'{e}t\'{e}s de rang 
 d'une partition de matrices  persym\'{e}triques.
 Nous montrons \'{e}galement que le nombre R de repr\'{e}sentations dans $ \mathbb{F}_{2}[T] $ de 0 comme 
 une somme de formes quadratiques associ\'{e}es aux sommes exponentielles  $\sum_{degY\leq  k-1}\sum_{degZ\leq s-1}E(tYZ) $
  est donn\'{e} par une int\'{e}grale   \'{e}tendue \`{a}  la boule unit\'{e} et est une combinaison lin\'{e}aire des $\Gamma _{i}.$
  Nous calculons alors explicitement le nombre R.
  Des r\'{e}sultats similaires sont \'{e}galement obtenus pour les $\mathbb{K} $- espaces vectoriels de dimension n+1.
  Nous terminons notre article en calculant explicitement le nombre de matrices de rang i de
   la forme  $\left[A\over B\right], $ o\`{u} A est persym\'{e}trique.
 \end{abstract}

\selectlanguage{english}

\begin{abstract}
 Let $\mathbb{K} $ be the field of Laurent Series  $ \mathbb{F}_{2}((T^{-1})). $\\
  We compute in particular   exponential sums in $\mathbb{K} $  of the form \\
   $  \sum_{degY\leq  k-1}\sum_{degZ\leq s-1}E(tYZ) $ where 
 t is in the unit interval of  $\mathbb{K},$  by showing that they only depend 
 on the rank  of some associated persymmetric matrices  with entries in  $\mathbb{F}_{2} $.
 ( A matrix\;$ [\alpha _{i,j}]  $ is  persymmetric if $ \alpha _{i,j} = \alpha _{r,s} $ \; for  \; i+j = r+s ). 
  Besides we establish rank properties of a partition of persymmetric matrices.
 We use these  results  to compute the number $ \Gamma_{i} $ of persymmetric matrices over   $ \mathbb{F}_{2} $ of
 rank i.
 We recover in this particular a general formula given by D. E. Daykin.
  Our proof is as indicated very different since it relies on rank properties 
 of a partition of persymmetric matrices. We also  prove that the number R of representations in  $ \mathbb{F}_{2}[T] $
  of 0 as a sum of some  quadratic forms associated to the exponential sums $\sum_{degY\leq  k-1}\sum_{degZ\leq s-1}E(tYZ) $
   is given by an integral over the unit interval, and is a linear combination of the $ \Gamma_{i}' s. $ We  then compute 
 explicitly the number R.
 Similar results are also obtained for n+1 dimensional $\mathbb{K} $ - vector spaces.We finish the paper by computing 
 explicitely the number of rank i matrices of the form $\left[A\over B\right], $ where A is  persymmetric.
 \end{abstract}

\newpage

\maketitle 
\newpage
\tableofcontents
\newpage
\section{Notations}
\label{sec 1}
\subsection{Analysis on   $\mathbb{K} $ }
\label{subsec 1.1}
We denote by $ \mathbb{F}_{2}\big(\big({\frac{1}{T}}\big) \big)
 = \mathbb{K} $ the completion
 of the field $\mathbb{F}_{2}(T), $  the field of  rational fonctions over the
 finite field\; $\mathbb{F}_{2}$,\; for the  infinity  valuation \;
 $ \mathfrak{v}=\mathfrak{v}_{\infty }$ \;defined by \;
 $ \mathfrak{v}\big(\frac{A}{B}\big) = degB -degA $ \;
 for each pair (A,B) of non-zero polynomials.
 Then every element non-zero t in
  $\mathbb{F}_{2}\big(\big({\frac{1}{T}}\big) \big) $
 can be expanded in a unique way in a convergent Laurent series
                              $  t = \sum_{j= -\infty }^{-\mathfrak{v}(t)}t_{j}T^j
                                 \; where\; t_{j}\in \mathbb{F}_{2}. $\\
  We associate to the infinity valuation\; $\mathfrak{v}= \mathfrak{v}_{\infty }$
   the absolute value \; $\vert \cdot \vert_{\infty} $\; defined by \;
  \begin{equation*}
  \vert t \vert_{\infty} =  \vert t \vert = 2^{-\mathfrak{v}(t)}. \\
\end{equation*}
    We denote  E the  Character of the additive locally compact group
$  \mathbb{F}_{2}\big(\big({\frac{1}{T}}\big) \big) $ defined by \\
\begin{equation*}
 E\big( \sum_{j= -\infty }^{-\mathfrak{v}(t)}t_{j}T^j\big)= \begin{cases}
 1 & \text{if      }   t_{-1}= 0, \\
  -1 & \text{if      }   t_{-1}= 1.
    \end{cases}
\end{equation*}
  We denote $\mathbb{P}$ the valuation ideal in $ \mathbb{K},$ also denoted the unit interval of  $\mathbb{K},$ i.e.
  the open ball of radius 1 about 0 or, alternatively, the set of all Laurent series 
   $$ \sum_{i\geq 1}\alpha _{i}T^{-i}\quad (\alpha _{i}\in  \mathbb{F}_{2} ) $$ and, for every rational
    integer j,  we denote by $\mathbb{P}_{j} $
     the  ideal $\left\{t \in \mathbb{K}|\; \mathfrak{v}(t) > j \right\}. $
     The sets\; $ \mathbb{P}_{j}$\; are compact subgroups  of the additive
     locally compact group \; $ \mathbb{K}. $\\
      All $ t \in \mathbb{F}_{2}\Big(\Big(\frac{1}{T}\Big)\Big) $ may be written in a unique way as
$ t = [t] + \left\{t\right\}, $ \;  $  [t] \in \mathbb{F}_{2}[T] ,
 \; \left\{t\right\}\in \mathbb{P}  ( =\mathbb{P}_{0}). $\\
 We denote by dt the Haar measure on \; $ \mathbb{K} $\; chosen so that \\
  $$ \int_{\mathbb{P}}dt = 1. $$\\
 \begin{defn}
\label{defn 1.1}We introduce the following definitions in  $ \mathbb{K}. $
\begin{itemize}
\item  A matrix\;$ D = [\alpha _{i,j}]  $ is said to be persymmetric if
$ \alpha _{i,j} = \alpha _{r,s} $ \; whenever \; i+j = r+s. \vspace{0.5 cm}
\item  We denote by  $ Persym _{n\times m}(\mathbb{F}_{2}) $  the set  of all 
 $ n \times m $  persymmetric matrices  over  $\mathbb{F}_{2}.$ \vspace{0.5 cm}
\item Let $ \Gamma _{i}^{s\times k} $ denote the number of persymmetric $ s\times k $ matrices  over $ \mathbb{F}_{2} $ of rank i.
We can assume without restriction that $ s\leq k,$  since obviously the transpose of a 
persymmetric matrix remains  persymmetric. \vspace{0.5 cm}
\item Set $ t = \sum_{i\geq 1}\alpha _{i}T^{-i}\in \mathbb{P} $  and let   $ (l,n,m) \in \mathbb{N^*}\times \mathbb{N^*}\times\mathbb{N^*}. $
We denote by $ D^{l}_{n\times m}(t) $   the following $ n\times m $ persymmetric matrix 
 $$ \left ( \begin{array} {cccccc}
\alpha _{l} & \alpha _{l+1} & \alpha _{l+2} &  \ldots & \alpha _{l+m-2}  &  \alpha _{l+m-1} \\
\alpha _{l+1 } & \alpha _{l+2} & \alpha _{l+3}&  \ldots  &  \alpha _{l+m-1} &  \alpha _{l+m} \\
\vdots & \vdots & \vdots   &  \ldots  & \vdots  & \vdots \\
\vdots & \vdots & \vdots    &  \ldots & \vdots & \vdots \\
\alpha _{l+n-2} & \alpha _{l+n-1} & \alpha _{l+n} & \ldots  &  \alpha _{l+n+m-4} &  \alpha _{l+n+m-3}  \\
\alpha _{l+n-1} & \alpha _{l+n} & \alpha _{l+n+1} & \ldots  &  \alpha _{l+n+m-3} &  \alpha _{l+n+m-2}  \\
\end{array}  \right). $$
If l =1 we denote $ D^{l}_{n\times m}(t) $  by  $ D_{n\times m}(t). $ \vspace{0.5 cm}
\item  We denote by ker D the Nullspace of  the matrix D and r(D) the rank of the matrix D.\vspace{0.5 cm}
\item Let $ h_{s,k}(t) = h(t)  $ be the quadratic exponential sum  defined  for rational  integers  $ s,k \geq 1 $ by
$$ t\in\mathbb{P} \longmapsto  \sum_{deg Y\leq k-1}\sum_{deg Z\leq s-1}E(tYZ) \in \mathbb{Z}. $$\vspace{0.5 cm}
\item  Let $ g_{s,k}(t) = g(t)  $ be the quadratic exponential sum  defined for rational  integers  $ s,k \geq 2 $ by
$$ t\in\mathbb{P}\longmapsto \sum_{deg Y= k-1}\sum_{deg Z=  s-1}E(tYZ) \in \mathbb{Z}. $$\vspace{0.5 cm}
 \item We denote by $ \mathbb{P}_{l-1}/\mathbb{P}_{l+n+m-2}$ a complete set
 of coset representatives of  $\mathbb{P}_{l+n+m-2}\; in\; \mathbb{P}_{l-1}. $
Obviously the sums  h(t) and g(t) are constant on cosets of $ \mathbb{P}_{k+s-1}. $\vspace{0.5 cm} 
 \item    Consider the following partition of the matrix $ D^{l}_{n\times m}(t), $
 obtained by drawing a horizontal line between the (n-1)th  and the nth row and drawing a vertical line between 
  the (m-1)th and the mth column in the matrix  $ D^{l}_{n\times m}(t). $
 $$ \left ( \begin{array} {ccccc|c}
\alpha _{l} & \alpha _{l+1} & \alpha _{l+2} &  \ldots & \alpha _{l+m-2}  &  \alpha _{l+m-1} \\
\alpha _{l+1 } & \alpha _{l+2} & \alpha _{l+3}&  \ldots  &  \alpha _{l+m-1} &  \alpha _{l+m} \\
\vdots & \vdots & \vdots   &  \ldots  & \vdots  & \vdots \\
\vdots & \vdots & \vdots    &  \ldots & \vdots & \vdots \\
\alpha _{l+n-2} & \alpha _{l+n-1} & \alpha _{l+n} & \ldots  &  \alpha _{l+n+m-4} &  \alpha _{l+n+m-3}  \\
\hline
\alpha _{l+n-1} & \alpha _{l+n} & \alpha _{l+n+1} & \ldots  &  \alpha _{l+n+m-3} &  \alpha _{l+n+m-2}  \\
\end{array}  \right). $$ 
$ Let  ( j_{1}, j_{2},   j_{3}, j_{4}) \in \mathbb{N}^{4}, $ we define \\
\begin{align*}
{}^{\#}\Big(\begin{array}{c | c}
           j_{1} & j_{2} \\
           \hline
           j_{3} & j_{4}
           \end{array} \Big)_{\mathbb{P}_{l-1}/\mathbb{P}_{l+n+m-2}}
            \end{align*}
    to be the cardinality of the following set
    $$\begin{array}{l}\{ t \in \mathbb{P}_{l-1}/\mathbb{P}_{l+n+m-2}
\mid r( D_{(n-1)\times( m-1)}^l(t) ) = j_{1} \; ,
r(D_{(n-1) \times m}^l(t)) = j_{2}  \\
 r( D_{n\times( m-1)}^l(t) ) = j_{3} \;  , r(D_{n \times m}^l(t)) = j_{4} \}. \end{array} $$\\
\end{itemize}
\end{defn}
\subsection{Analysis on the (n+1) - dimensional   $\mathbb{K}$- vectorspace}
\label{subsec 1.2}
  Let $ \mathbb{K}^{n+1} $ be the (n+1)- dimensional vector space  over $ \mathbb{K}. $
  Let $ (t,\eta_{1},\eta_{2},\ldots,\eta_{n} )\in \mathbb{K}^{n+1}  $ 
   and  $ \vert (t,\eta_{1},\eta_{2},\ldots,\eta_{n} )\vert =
  sup \left\{\vert t \vert ,\vert \eta_{1} \vert, \vert \eta_{2} \vert,\ldots,\vert \eta_{n} \vert\right\}
   =  2^{-inf(\mathfrak{v}(t),\mathfrak{v}(\eta_{1} ),\mathfrak{v}(\eta_{2} ),\ldots,\mathfrak{v}(\eta_{n} ))}. $\\
   
 It is easy to see that $(t,\eta_{1},\eta_{2},\ldots,\eta_{n} )\longrightarrow \vert (t,\eta_{1},\eta_{2},\ldots,\eta_{n} )\vert $
  is an ultrametric valuation on $ \mathbb{K}^{n+1}, $ that is \\
  
  $(t,\eta_{1},\eta_{2},\ldots,\eta_{n} )\longrightarrow \vert (t,\eta_{1},\eta_{2},\ldots,\eta_{n} )\vert $ is a norm and
  
  $  \vert((t,\eta_{1},\eta_{2},\ldots,\eta_{n} ) + (t',\eta_{1}',\eta_{2}',\ldots,\eta_{n}' ))\vert \leq max\left\{\vert  (t,\eta_{1},\eta_{2},\ldots,\eta_{n} )\vert ,\vert (t',\eta_{1}',\eta_{2}',\ldots,\eta_{n}' )\vert \right\}. $ \\
  
   We denote by  $ d((t,\eta_{1},\eta_{2},\ldots,\eta_{n} ) =dtd\eta_{1}d\eta_{2} \cdots d\eta_{n} $ the  Haar measure on   $ \mathbb{K}^{n+1} $
  chosen so that the measure on the unit interval of $ \mathbb{K}^{n+1}  $ is equal to one, that is 
  
 $$ \int_{\mathbb{P}^{n+1}} d((t,\eta_{1},\eta_{2},\ldots,\eta_{n} ) =
   \int_{\mathbb{P}}dt\int_{\mathbb{P}}d\eta_{1}\int_{\mathbb{P}}d\eta_{2}\cdots \int_{\mathbb{P}}d\eta_{n} = 1\cdot 1\cdot1\cdots 1 =1. $$
 
Let  $$ (t,\eta_{1},\eta_{2},\ldots,\eta_{n} )
 =  \big( \sum_{i= -\infty }^{-\mathfrak{v}(t)}t_{i}T^i, \sum_{i = -\infty }^{-\mathfrak{v}(\eta_{1} )}\eta _{1i}T^i,\sum_{i = -\infty }^{-\mathfrak{v}(\eta_{2} )}\eta _{2i}T^i,\ldots,\sum_{i = -\infty }^{-\mathfrak{v}(\eta_{n} )}\eta _{ni}T^i   \big) \in  \mathbb{K}^{n+1}, $$ 
 we denote $\chi $  the  Character on  $(\mathbb{K}^{n+1}, +) $ defined by \\
 
\begin{align*}
& \chi\big( \sum_{i= -\infty }^{-\mathfrak{v}(t)}t_{i}T^i ,
\sum_{ i = -\infty }^{-\mathfrak{v}(\eta_{1} )}\eta _{1i}T^i,  \sum_{ i = -\infty }^{-\mathfrak{v}(\eta _{2})}\eta _{2i}T^i,\ldots,\sum_{ i = -\infty }^{-\mathfrak{v}(\eta_{n} )}\eta _{ni}T^i  \big) \\
& =  E\big( \sum_{i = -\infty }^{-\mathfrak{v}(t)}t_{i}T^i\big)\cdot
  E\big( \sum_{i = -\infty }^{-\mathfrak{v}(\eta_{1} )}\eta _{1,i}T^i\big)\cdot E\big( \sum_{i = -\infty }^{-\mathfrak{v}(\eta_{2} )}\eta _{2,i}T^i\big)\cdots E\big( \sum_{i = -\infty }^{-\mathfrak{v}(\eta_{n} )}\eta _{n,i}T^i\big) \\
& = \begin{cases}
 1 & \text{if      }   t_{-1} +  \eta _{1,-1}+ \eta _{2,-1}+\ldots + \eta _{n,-1}= 0, \\
  -1 & \text{if      }   t_{-1} + \eta _{1,-1}+\eta _{2,-1}+\ldots + \eta _{n,-1}= 1.
    \end{cases}
\end{align*}
\begin{defn}
\label{defn 1.2}We introduce the following definitions in the (n+1) - dimensional  $ \mathbb{K} $- vectorspace.
\begin{itemize}
\item  We denote by $ \mathbb{P}/\mathbb{P}_{k+m} \times ( \mathbb{P}/\mathbb{P}_{k})^{n} $
  a complete set of coset representatives of
$\mathbb{P}_{k+m} \times \mathbb{P}_{k}^{n}\; in \;  \mathbb{P}^{n+1}. $ \vspace{0.5 cm}
\item  $$ Set\quad
   (t,\eta_{1},\eta_{2},\ldots,\eta_{n} )
 =  \big( \sum_{i\geq 1}\alpha _{i}T^i, \sum_{i \geq 1}\beta  _{1i}T^i,\sum_{i \geq 1}\beta  _{2i}T^i,\ldots,\sum_{i \geq 1}\beta _{ni}T^i   \big) \in  \mathbb{P}^{n+1}. $$ 
 We denote by $ D^{\Big[\stackrel{n}{1+m}\Big] \times k } (t,\eta_{1},\eta_{2},\ldots,\eta_{n} ) $
    the following $(1+n+m)\times k $ matrix, where the submatrix obtained by 
  delating the n last rows form a $ (1+m)\times k $ persymmetric matrix, and the n last rows form a $n\times k $ matrix
  over the finite field  $\mathbb{F}_{2} $ 
    $$  \left ( \begin{array} {cccccc}
\alpha _{1} & \alpha _{2} & \alpha _{3} &  \ldots & \alpha _{k-1}  &  \alpha _{k} \\
\alpha _{2 } & \alpha _{3} & \alpha _{4}&  \ldots  &  \alpha _{k} &  \alpha _{k+1} \\
\vdots & \vdots & \vdots   &  \ldots  & \vdots  &  \vdots \\
\vdots & \vdots & \vdots    &  \ldots & \vdots  &  \vdots \\
\alpha _{1+m} & \alpha _{2+m} & \alpha _{3+m} & \ldots  &  \alpha _{k+m-1} &  \alpha _{k+m}  \\
\hline
\beta  _{11} & \beta  _{12} & \beta  _{13} & \ldots  &  \beta_{1 k-1} &  \beta _{ 1k}  \\
\beta  _{21} & \beta  _{22} & \beta  _{23} & \ldots  &  \beta_{2 k-1} &  \beta _{2k}\\
\vdots & \vdots & \vdots   &  \ldots  & \vdots  &  \vdots \\
\vdots & \vdots & \vdots    &  \ldots & \vdots  &  \vdots \\
\beta  _{n1} & \beta  _{n2} & \beta  _{n3} & \ldots  &  \beta_{n k-1} &  \beta _{nk}
 \end{array}  \right). $$\vspace{0.5 cm}
  \item  Let $ g_{m,k}(t,\eta ) = g(t,\eta ) $ be the exponential sum defined for rational  integers  $ m\geq  0, 
k\geq 1 $ by
$$ (t,\eta ) \in  \mathbb{P}\times \mathbb{P}\longmapsto  
  \sum_{deg Y\leq k-1}\sum_{deg Z\leq m}E(tYZ)\sum_{deg U =0}E(\eta YU) \in \mathbb{Z}.  $$\vspace{0.5 cm}
\item Let $ f_{m,k}(t,\eta ) = f(t,\eta ) $ be the exponential sum defined for rational  integers  $ m\geq  0, 
k\geq 1 $ by 
$$ (t,\eta ) \in  \mathbb{P}\times \mathbb{P}\longmapsto 
  \sum_{deg Y\leq k-1}\sum_{deg Z\leq m}E(tYZ)\sum_{deg U \leq 0}E(\eta YU) \in \mathbb{Z}.  $$\vspace{0.5 cm}
\item  \[\begin{array}{lcl}
 \sigma _{i,i}^{\big[\stackrel{1}{1+m}\big]\times k} &  = & Card
  \left\{(t,\eta )\in \mathbb{P}/\mathbb{P}_{k+m}\times \mathbb{P}/\mathbb{P}_{k}
\mid  r(D_{(1+m)\times k}(t)) = r(D^{\big[\stackrel{1}{1+m}\big] \times k }(t,\eta )) = i
 \right\},\\
  \sigma _{i-1,i}^{\big[\stackrel{1}{1+m}\big]\times k} &  = & Card
  \left\{(t,\eta )\in \mathbb{P}/\mathbb{P}_{k+m}\times \mathbb{P}/\mathbb{P}_{k}
\mid  r(D_{(1+m)\times k}(t)) = i-1, \;  r(D^{\big[\stackrel{1}{1+m}\big] \times k }(t,\eta )) = i
 \right\},\\
  \Gamma _{i}^{\Big[\substack{n \\ 1+m }\Big] \times k} & =  & Card
  \left\{(t,\eta_{1}, \eta_{2},\ldots,\eta_{n})\in \mathbb{P}/\mathbb{P}_{k+m}\times (\mathbb{P}/\mathbb{P}_{k})^{n}
\mid   r(D^{\big[\stackrel{n}{1+m}\big] \times k })(t,\eta_{1}, \eta_{2},\ldots,\eta_{n}) = i
 \right\}\end{array}\]
\end{itemize}
\end{defn} 

 \section{Introduction (refer to Section \ref{sec 1} and Section \ref{sec 3})}
\label{sec 2}
The rational function field  $\mathbb{F}_{2}(T) $ is completed with respect to an appropriate valuation 
to a field $ \mathbb{K} $ (i.e. the field of Laurent Series). The unit interval of $\mathbb{K},$ that is, the open ball
of radius 1 about 0, is a compact additive group. We shall use the Haar integral on this group.\\
The paper is based on the proof of the following  formula  of the number $ \Gamma_{i} $
 of $ s\times k $ persymmetric matrices over $\mathbb{F}_{2} $ of rank i, that is 
 \begin{equation}
 \label{eq 2.1}
\Gamma_{i}^{s\times k} = \left\{\begin{array}{ccc}
             1 & if & i= 0, \\
             3\cdot 2^{2(i-1)} & if & 1\leq i \leq s-1, \\
              2^{k+s-1} - 2^{2s-2} & if & i= s\; (\leq k ).
             \end{array}\right.\
\end{equation}
We remark that David E. Daykin [1] has already proved this result over any finite field $\mathbb{F}$
with the number 2 in the formula replaced by $ |\mathbb{F}|, $   and the number 3
replaced by $ |\mathbb{F}|^2 -1. $   Our proof is very different since it relies on rank properties 
 of a partition of persymmetric matrices, and is proper to the finite field with two elements. Besides, in order  to establish
the formula \eqref{eq 2.1}, we obtain several results  concerning exponential sums in  $ \mathbb{K}$ and persymmetric matrices,
 which seem to have an  intrinsic interest. \\
\textbf{ \underline {The paper is organized as follows:}}\vspace{0.2 cm}\\
 \textbf{ In Section 1,} are introduced main notations and definitions. \\
  \textbf{ In Section 3,} we state the main theorems  in the  (n+1) dimensional $ \mathbb{K} $-vector spaces.\\
   \textbf{ In Section 4,}    we establish  in particular  results  on  exponential sums in $ \mathbb{P}$ of the form h(t), g(t) and we show
    that they only depend on rank properties of some corresponding  persymmetric matrices with entries in $ \mathbb{F}_{2}.$ \\
The proof of these results is based on the following identity \\
\begin{equation}
\label{eq 2.2}
   \sum_{deg Z\leq s-1} E(tYZ)   =  \left\{ \begin{array}{ccc}
                         2^{s}    &  if  & Y\in \ker D_{s\times k}(t) \\
                          0  &       &  otherwise,
                          \end{array}\right.  
\end{equation}
 where  $ t \in\mathbb{P}\; and \; \deg Y\leq k-1 $.\\
By \eqref{eq 2.2} we obtain \\
\begin{equation}
\label{eq 2.3}
 Let \; t\in\mathbb{P},\quad  then \; h_{s,k}(t) = h(t)= \sum_{deg Y\leq k-1}\sum_{deg Z\leq s-1}E(tYZ) \quad is \quad given \quad by \quad
2^{k+s-r(D_{s\times k}(t))}.
\end{equation}
From \eqref{eq 2.3}, observing that  h(t)  is  constant on cosets of $ \mathbb{P}_{k+s-1} $, we  prove that the number R of representations
 in $ \mathbb{F}_{2}[T],$ of 0 as a sum of 
quadratic forms  is given by an integral over the unit interval of  $\mathbb{K},$
and is a linear combination of the  $ \Gamma_{i}'s. $ \\
More precisely  we prove 
\begin{equation}
\label{eq 2.4}
  R =  \int_{\mathbb{P}}h^{q}(t) dt 
 = 2^{(q-1)(k+s) +1}\sum_{i=0}^{s}\Gamma_{i}2^{-qi}.
\end{equation}
By calculating the square of g(t) and using \eqref{eq 2.3} we obtain  \\
\begin{equation}
\label{eq 2.5}
 Let \quad  t\in\mathbb{P},\quad then \quad g_{s,k}(t) = g(t) = \sum_{deg Y= k-1}\sum_{deg Z = s-1}E(tYZ) \quad  is \quad given \quad by 
 \end{equation}
 \begin{equation*}
 \begin{cases}
 2^{s+k-j-2} &\text{if  }  r(D_{(s-1) \times (k-1) }(t)) = r(D_{s \times ( k-1)}(t)) = r(D_{(s-1) \times k }(t)) = r(D_{s \times k}(t)) = j, \\
  -  2^{s+k-j-2}  &\text{if  }  r(D_{(s-1) \times (k-1) }(t)) = r(D_{s \times ( k-1)}(t)) = r(D_{(s-1) \times k }(t)) =  j \;  and \;r(D_{s \times k}(t)) = j+1, \\
   0     &  \text{if  }    otherwise.
      \end{cases}
\end{equation*}
\textbf{In Section 5,} we study  rank properties of a partition  of persymmetric matrices
 by  considering  the following partition of the matrix $ D_{s\times k}(t) $
  $$ \left ( \begin{array} {ccccc|c}
\alpha _{1} & \alpha _{2} & \alpha _{3} &  \ldots & \alpha _{k-1}  &  \alpha _{k} \\
\alpha _{2 } & \alpha _{3} & \alpha _{4}&  \ldots  &  \alpha _{k} &  \alpha _{k+1} \\
\vdots & \vdots & \vdots   &  \vdots  & \vdots  & \vdots \\
\vdots & \vdots & \vdots    &  \vdots & \vdots & \vdots \\
\alpha _{s-1} & \alpha _{s} & \alpha _{s+1} & \ldots  &  \alpha _{k+s-3} &  \alpha _{k+s-2}  \\
\hline
\alpha _{s} & \alpha _{s+1} & \alpha _{s+2} & \ldots  &  \alpha _{k+s-2} &  \alpha _{k+s-1}  \\
\end{array}  \right). $$ \\
We prove \\

 \begin{align}  If\quad  0\leq j\leq s-1, \quad then \quad
{}^{\#}\Big(\begin{array}{c | c}
           j & j \\
           \hline
           j  &  j
           \end{array} \Big)_{\mathbb{P}/\mathbb{P}_{k +s -1}}  =
           {}^{\#}\Big(\begin{array}{c | c}
           j & j \\
           \hline
           j  &  j +1
           \end{array} \Big)_{\mathbb{P}/\mathbb{P}_{k +s -1}}.\label{eq 2.6}
            \end{align}\\
            This result is deduced from the fact that 
              \[\int_{\mathbb{P}}g^{2q+1}(t)dt = 0. \] 
              
 \begin{equation}
 \label{eq 2.7}
  If\quad 0\leq j\leq s-2, \quad then \quad
{}^{\#}\Big(\begin{array}{c | c}
           j & j+1 \\
           \hline
           j +1 &  j +1
           \end{array} \Big)_{\mathbb{P}/\mathbb{P}_{k +s -1}}  = 0 . 
            \end{equation}
            Actually, this result holds  for any matrix.
            
  \begin{align}
  If\quad 0\leq j\leq s-2,\quad then \quad 
{}^{\#}\Big(\begin{array}{c | c}
           j & j +1 \\
           \hline
           j  &  j + 1
           \end{array} \Big)_{\mathbb{P}/\mathbb{P}_{k +s -1}}  =
           {}^{\#}\Big(\begin{array}{c | c}
           j & j \\
           \hline
           j +1 &  j +1
           \end{array} \Big)_{\mathbb{P}/\mathbb{P}_{k +s -1}} = 0.\label{eq 2.8}
            \end{align}
           \vspace{0.1 cm} 
            These equalities are consequence of rank properties of submatrices of persymmetric $ s\times k $ matrices.\\
            
          By combining \eqref{eq 2.6},\eqref{eq 2.7} and \eqref{eq 2.8}  we obtain\\
          \begin{equation}
          \label{eq 2.9}
  2\cdot\Gamma_{i}^{s \times (k-1)} = 2\cdot\Gamma_{i}^{(s-1) \times k} 
 = 2\cdot {}^{\#}\Big (\begin{array}{c| c}
        i & i \\
       \hline
       i & i
\end{array}\Big)_{ \mathbb{P}/\mathbb{P}_{k+s-1}}
+{}^{\#}\Big (\begin{array}{c| c}
        i-1 & i \\
       \hline
       i & i+1
\end{array}\Big)_{ \mathbb{P}/\mathbb{P}_{k+s-1}}.
\end{equation}

From \eqref{eq 2.9} we get \\

 \begin{equation}
        \label{eq 2.10}
         If \quad  0\leq i \leq s-2, \quad  then \quad    
    \Gamma_{i}^{s \times (k-1)} =   \Gamma_{i}^{(s-1) \times k}. 
       \end{equation}
   \textbf{ In Section 6,}  
  we compute exponential sums  in  $\mathbb{P}^{2} $ associated to a matrix of the form  $\big[{A\over b_{-}}\big],$ 
 where A is a $ (m+1)\times k $ persymmetic matrix over  $ \mathbb{F}_{2} $ and $ b_{-} $ a $1\times k $ matrix with entries in  $ \mathbb{F}_{2}.$
 
 In particular we prove that the number of rank i  matrices of the form  \\
   
  \begin{center}
    \(\left( \begin{array} {cccccc}
\alpha _{1} & \alpha _{2} & \alpha _{3} &  \ldots & \alpha _{k-1}  &  \alpha _{k} \\
\alpha _{2 } & \alpha _{3} & \alpha _{4}&  \ldots  &  \alpha _{k} &  \alpha _{k+1} \\
\vdots & \vdots & \vdots   &  \ldots  & \vdots  &  \vdots \\
\vdots & \vdots & \vdots    &  \ldots & \vdots  &  \vdots \\
\alpha _{1+m} & \alpha _{2+m} & \alpha _{3+m} & \ldots  &  \alpha _{k+m-1} &  \alpha _{k+m}  \\
\hline
\beta  _{1} & \beta  _{2} & \beta  _{3} & \ldots  &  \beta_{k-1} &  \beta _{k}  
 \end{array}\right)\)\quad 
 \end{center}  
  is a linear combination of the $ \Gamma _{j}^{(1+m)\times k} $ for $ j \in\{i-1, i \} $
 ( recall that  $ \Gamma _{j}^{(1+m)\times k} $ denotes  the number of persymmetric $(1+m)\times k $ matrices of rank j ).\\
  In fact we prove that \\
  \begin{equation}
\label{eq 2.11}
   \Gamma _{i}^{\Big[\substack{1 \\ 1+m }\Big] \times k} = ( 2^{k} -  2^{i-1})\cdot\Gamma _{i-1}^{(1+m)\times k} 
 +  2^{i}\Gamma _{i}^{(1+m)\times k} \quad for \quad all \quad 1\leq i\leq inf(k,2+m). 
\end{equation}
 The proof of \eqref{eq 2.11} is based on the following identity \\
 \begin{equation} 
 \label{eq 2.12}
E(\eta Y) = \begin{cases}
1 & \text{if        }                      \sum_{j=1}^{k-1}\beta _{j}\delta _{j} = 0 \\
-1 & \text{if        }                       \sum_{j=1}^{k-1}\beta _{j}\delta _{j} = 1
\end{cases}
\end{equation}
 where,  $$  \eta  = \sum_{j\geq 1}\beta _{j}T^{-j} \in \mathbb{P}\quad
  and \quad Y = \sum_{i=1}^{k-1}\delta _{i}T^i\in \mathbb{F}_{2}[T] ,\quad degY\leq k-1 . $$ \\
  
  We proceed as follows : \\
  
   By \eqref{eq 2.12} we obtain \\
   
  \begin{equation}
  \label{eq 2.13}
 g(t,\eta ) = \begin{cases}
 2^{k+m+1-  r(D_{(1+m)\times k}(t)) }  & \text{if }
  r(D_{(1+m)\times k}(t)) = r(D^{\big[\stackrel{1}{1+m}\big] \times k }(t,\eta ) ), \\
     0  & \text{otherwise },
    \end{cases}
\end{equation}
$$  where \quad    g(t,\eta ) = \sum_{deg Y\leq k-1}\sum_{deg Z\leq m}E(tYZ)\sum_{deg U = 0}E(\eta YU) .$$
\begin{align*}
& 
\end{align*}

 By squaring $ g(t,\eta ) $ we obtain  $ g^2(t,\eta ) = g(t,\eta )\cdot h(t) $
   where $ h(t)  = \sum_{deg \leq k-1}\sum_{deg Z \leq m}E(tYZ) $.\\
   By recurrence on q we get \\
   \begin{equation}
\label{eq 2.14}
  g^{q}(t,\eta ) = g(t,\eta )\cdot h^{q-1}(t).
\end{equation}
\begin{align*}
& 
\end{align*}
 By integrating $ g^q(t,\eta )  $ over the unit interval of $ \mathbb{K}^{2},$ observing that $ g(t,\eta ) $ is constant on 
  cosets of  $ \mathbb{P}_{k+m}\times\mathbb{P}_{k} $ and using \eqref{eq 2.13},  we obtain\\
  \begin{equation}
\label{eq 2.15}
   \int_{\mathbb{P}}\int_{\mathbb{P}}g^q(t,\eta ) dt d\eta  =
2^{q(k+m +1) -2k - m}\sum_{i = 0}^{inf(k,1+m)}
 \sigma _{i,i}^{\big[\stackrel{1}{1+m}\big]\times k}2^{-iq}. 
\end{equation}
\begin{align*}
& 
\end{align*}
 Integrating $  g(t,\eta )\cdot h^{q-1}(t)  $ over the unit interval of $ \mathbb{K}^{2}, $  we get by Fubini's theorem \\
 \begin{equation}
 \label{eq 2.16}
   \int_{\mathbb{P}}\int_{\mathbb{P}}g(t,\eta )h^{q-1}(t) dt d\eta    =
2^{q(k+m +1) -2k - m}\sum_{i = 0}^{inf(k,1+m)}2^{i}\Gamma _{i}^{(1+m)\times k}
2^{-iq}.  
\end{equation}
\begin{align*}
& 
\end{align*}
  Comparing \eqref{eq 2.15} and \eqref{eq 2.16},  using \eqref{eq 2.14}, we get \\
 \begin{equation}
  \label{eq 2.17}
  \sigma _{i,i}^{\big[\stackrel{1}{1+m}\big]\times k} =
 2^{i}\Gamma _{i}^{(1+m)\times k} \quad for\; all \; 0\leq i\leq inf(k, 1+m). 
  \end{equation}
\begin{align*}
& 
\end{align*}
   We obviously have \\
   \begin{equation}
   \label{eq 2.18}
    \Gamma _{i}^{\Big[\substack{1 \\ 1+m }\Big] \times k}
=    \sigma _{i,i}^{\big[\stackrel{1}{1+m}\big]\times k}+ \sigma _{i-1,i}^{\big[\stackrel{1}{1+m}\big]\times k}.    
\end{equation}
\begin{align*}
& 
\end{align*}
  From \eqref{eq 2.17} and \eqref{eq 2.18} we  easily get \eqref{eq 2.11}.
  \begin{align*}
& 
\end{align*}
   Besides we prove, similarly to \eqref{eq 2.3} and \eqref{eq 2.4} that \\
 \begin{align}
&  f(t,\eta ) = 2^{k+m +2 -r(D^{\big[\stackrel{1}{1+m}\big] \times k }(t,\eta )) }\label{eq 2.19}\\
 & \text{and}\quad  \int_{\mathbb{P}}\int_{\mathbb{P}}f^q(t,\eta ) dt d\eta  =
2^{q(k+m +2) -2k - m}\sum_{i = 0}^{inf(k,2+m)} \Gamma _{i}^{\Big[\substack{1 \\ 1+m }\Big] \times k}2^{-iq} \label{eq 2.20}\\
& \text{where} \quad  f(t,\eta ) =  \sum_{deg Y\leq k-1}\sum_{deg Z\leq m}E(tYZ)\sum_{deg U \leq 0}E(\eta YU). \nonumber \\
& \nonumber
 \end{align}
 \textbf{ In Section 7,}   we prove  formula \eqref{eq 2.1}  by induction on $ s \;(\leq k) $ and by using \\
    \begin{equation*}
\eqref{eq 2.9}\quad  \Gamma_{i}^{s \times (k-1)} =   \Gamma_{i}^{(s-1) \times k} \quad for \quad  0\leq i\leq  s-2 
\end{equation*}
  and the two following equations \\
   \begin{align*}
 \sum_{i=0}^{s} \Gamma_{i}^{s \times k}& =  2^{k+s-1}, \\
 \sum_{i=0}^{s}\Gamma_{i}^{s \times k}\cdot2^{-i} &  =   2^{k-1} +2^{s-1} -2^{-1}. \\
 & 
 \end{align*}
    \textbf{ In Section 8} we prove Theorem \ref{thm 3.3} by showing that the solutions of the system  
    \begin{equation}
    \label{eq 2.21}
     \Gamma _{i}^{s\times k}=
   \sum_{(j_{1}, j_{2},j_{3})\in \{i-2, i-1, i\}^{3}} {}^{\#}\Big(\begin{array}{c | c}
           j_{1} & j_{2} \\
           \hline
             j_{3}  &  i
           \end{array} \Big)_{\mathbb{P}/\mathbb{P}_{k +s -1}}  = \left\{\begin{array}{ccc}
             1 & if & i= 0, \\
             3\cdot 2^{2(i-1)} & if & 1\leq i \leq s-1, \\
              2^{k+s-1} - 2^{2s-2} & if & i= s\; (\leq k),
             \end{array}\right.\
 \end{equation}
  are given by 
  \begin{equation}
  \label{eq 2.22}
{}^{\#}\Big(\begin{array}{c | c}
           j_{1} & j_{2} \\
           \hline
           j_{3} & j_{4}
           \end{array} \Big)_{\mathbb{P}/\mathbb{P}_{k+s-1}} =
  \begin{cases}
  1 & \text{if  }  j_{1}= j_{2}= j_{3}= j_{4} = 0, \\
   2^{2j-1} & \text{if  }  j_{1}=j_{2}=j_{3}=j \; ,j_{4}\in\left\{j,j+1\right\} \;,1\leq j\leq s-1,\\
    2^{2j-3} & \text {if  }  j_{1}=  j-2 \;, j_{2}= j_{3}= j-1\;,  j_{4}= j \;, 2\leq j\leq s, \\
    2^{k+s-1}-2^{2s-1} & \text{if   }  j_{1}= j_{2}= s-1 ,  j_{3}= j_{4}= s \\
     0 & \text{otherwise}.
   \end{cases}
   \end{equation} 
 \textbf{ In Subsection 9.1},\quad we show that Theorem \ref{thm 3.4} follows from \eqref{eq 2.3} and \eqref{eq 2.4}.\\

\textbf{ In Subsection 9.2},\quad  we show that Theorem \ref{thm 3.5} follows from \eqref{eq 2.5} and \eqref{eq 2.10}.\\

\textbf{ In Subsection 9.3},\quad  we show that Theorem \ref{thm 3.6} follows from \eqref{eq 2.13} and \eqref{eq 2.15}.\\

\textbf{ In Subsection 9.4},\quad  we show that Theorem \ref{thm 3.7} follows from \eqref{eq 2.19} and \eqref{eq 2.20}.\\

\textbf{ In Subsection 9.5},\quad  we show that Theorem \ref{thm 3.8} follows from \eqref{eq 2.11} and \eqref{eq 2.1} with 
$s\rightarrow 1+m $. \\
 
     \textbf{ In Subsection 10.1 },
   we show that  $ \Gamma _{i}^{\Big[\substack{n \\ 1+m }\Big] \times k}$ can be written as  a linear
   combination of the  $\Gamma _{i-j}^{(1+m)\times k}$ for $ j = 0,1,\ldots,n, $
    where $ \Gamma _{i}^{\Big[\substack{n \\ 1+m }\Big] \times k}$ denotes the number of 
   rank i matrices of the form $ [{A\over B}] $ such that A is  a $(1+m)\times k$  persymmetric matrix
   and B is a $ n\times k $ matrix with entries in  $ \mathbb{F}_{2} $
  and  $\Gamma _{i-j}^{(1+m)\times k}$
  the number of $(1+m)\times k$ persymmetric matrices  of rank i-j.\\
   Precisely we show by induction on n \\
 
   \begin{equation}
   \label{eq 2.23}
   \Gamma _{i}^{\Big[\substack{n \\ 1+m }\Big] \times k}=
 \sum_{j= 0}^{n}\Big[2^{(n-j)\cdot(i-j)} a_{j}^{(n)}\prod_{l=1}^{j}(2^{k}- 2^{i-l})\Big]\cdot
 \Gamma _{i-j}^{(1+m)\times k} \quad for \quad 0\leq i\leq inf(k,n+m +1)
 \end{equation} where $ a_{j}^{(n)} $ satisfies the linear recurrence relation 
 \begin{equation*}
  a_{j}^{(n)} = 2^{j}\cdot a_{j}^{(n-1)} + a_{j-1}^{(n-1)},\quad n = 2,3,4,\ldots       \quad for\quad 1\leq j\leq n-1. 
\end{equation*}
 \textbf{ In Subsection 10.2 },  the explicit value of $ a_{j}^{(n)} $ in \eqref{eq 2.23}, is obtained as follows :\\
  Setting m = 0 in \eqref{eq 2.23} we have:
 \begin{equation}
\label{eq 2.24}
  \Gamma _{i}^{\Big[\substack{n \\ 1}\Big] \times k} = 
    \sum_{j= 0}^{n}2^{(n-j)\cdot(i-j)} a_{j}^{(n)}\prod_{l=1}^{j}(2^{k}- 2^{i-l})
 \Gamma _{i-j}^{ 1 \times k} \quad for \quad 0\leq i\leq inf(k,n+1). 
\end{equation}
 Recalling that  $ \Gamma _{i}^{\Big[\substack{n \\ 1 }\Big] \times k} $ denotes the number of
   $ (n+1)\times k $  matrices with entries in $ \mathbb{F}_{2} $ of rank i,we get  by combining \eqref{eq 2.24} and [3, George Landsberg]
  \begin{equation*}
 \Gamma _{i}^{\Big[\substack{n \\ 1 }\Big] \times k}=
 \prod_{l = 0}^{i-1}{ (2^{n+1} -2^{l})(2^{k}-2^{l}) \over (2^{i}- 2^{l})} 
  =  \sum_{j= i-1}^{i}2^{(n-j)\cdot(i-j)} a_{j}^{(n)}\prod_{l=1}^{j}(2^{k}- 2^{i-l}) 
   \Gamma _{i-j}^{1 \times k} 
\end{equation*}
   from which, we deduce 
     \begin{equation*}
  a_{i}^{(n)} +  2^{n-(i-1)}\cdot a_{i-1}^{(n)}=  \prod_{l = 0}^{i-1}{ 2^{n+1} -2^{l} \over 2^{i}- 2^{l}}. 
 \end{equation*} \\
 The formula for $ a_{i}^{(n)} $ given in Lemma \ref{lem 10.2} then follows.\\
  \textbf{ In Subsection 10.3 }, we obtain Theorem \ref{thm 3.9} from \eqref{eq 2.23} and the explicit value of  $ a_{j}^{(n)}. $ \\
  \textbf{ In Subsection 10.4 }, we get Corollary \ref{cor 3.10} by computing  $ a_{j}^{(n)} $ for $0 \leq j\leq n $,\quad $1\leq n\leq 5$.\\
   \textbf{ In Subsection 10.5 }, we generalize \eqref{eq 2.19} and \eqref{eq 2.20}.

 \section{Statement of results}
\label{sec 3}

 \begin{thm}
\label{thm 3.1}
     The number $ \Gamma_{i}^{s\times k} $ of persymmetric
      $ s\times k $ matrices over $\mathbb{F}_{2}$  of rank i  is given by  \\
 \[\Gamma_{i}^{s\times k} = \left\{\begin{array}{ccc}
             1 & if & i= 0, \\
             3\cdot 2^{2(i-1)} & if & 1\leq i \leq s-1, \\
          2^{k+s-1} - 2^{2s-2} & if &  i = s \; (s\leq k).
             \end{array}\right.\] \\
\end{thm}

\begin{remark}
\label{remark 3.2}
David E. Daykin has already proved this result over any finite field $\mathbb{F}$
with the number 2 in the formula replaced by$ |\mathbb{F}|,$   and the number 3
replaced by $ |\mathbb{F}|^2 -1 $, see [1].\;
Our proof is different and proper to the finite field with two elements.
\end{remark}
\begin{thm}
\label{thm 3.3}
 Let $ (j_{1},j_{2},j_{3},j_{4})\in \mathbb{N}^4 ,$  then
\begin{equation*}
{}^{\#}\Big(\begin{array}{c | c}
           j_{1} & j_{2} \\
           \hline
           j_{3} & j_{4}
           \end{array} \Big)_{\mathbb{P}/\mathbb{P}_{k+s-1}} =
  \begin{cases}
  1 & \text{if  }  j_{1}= j_{2}= j_{3}= j_{4} = 0, \\
   2^{2j-1} & \text{if  }  j_{1}=j_{2}=j_{3}=j \; ,j_{4}\in\left\{j,j+1\right\} \;,1\leq j\leq s-1,\\
    2^{2j-3} & \text {if  }  j_{1}=  j-2 \;, j_{2}= j_{3}= j-1\;,  j_{4}= j \;, 2\leq j\leq s, \\
    2^{k+s-1}-2^{2s-1} & \text{if   }  j_{1}= j_{2}= s-1 ,  j_{3}= j_{4}= s, \\
     0 & \text{otherwise},
   \end{cases}
   \end{equation*} \\
 where
  ${}^{\#}\Big(\begin{array}{c | c}
           j_{1} & j_{2} \\
           \hline
           j_{3} & j_{4}
           \end{array} \Big)_{\mathbb{P}/\mathbb{P}_{k+s -1}} $
          was  introduced  in  Definition \ref{defn 1.1}.
\end{thm}
  \begin{thm}
\label{thm 3.4}
 Let $ h_{s,k}(t) = h(t)  $ be the quadratic exponential sum in $ \mathbb{P} $ defined by
$$ t\in\mathbb{P}\longmapsto  \sum_{deg Y\leq k-1}\sum_{deg Z\leq s-1}E(tYZ) \in \mathbb{Z}. $$\vspace{0.5 cm}
 Then
$$ h(t) = 2^{k+s - r(D_{s\times k}(t))} $$ and \\
   $$   \int_{\mathbb{P}}h^{q}(t) dt = 2^{(q-1)(k+s) +1}\sum_{i=0}^{s}\Gamma_{i}^{s\times k}2^{-qi}.   $$
Let R denote the number of solutions
 $(Y_1,Z_1, \ldots,Y_q,Z_q) $  of the polynomial equation
                        $$ Y_1Z_1 +  Y_2Z_2 + \ldots + Y_qZ_q = 0  $$ \\
  satisfying the degree conditions \\
                   $$  degY_i \leq k-1 , \quad degZ_i \leq s-1 \quad for  1\leq i \leq q. $$ \\
 Then $$ R =  \int_{\mathbb{P}}h^{q}(t) dt $$    
     \end{thm}
\begin{thm}
\label{thm 3.5}
 Let $ g_{s,k}(t) = g(t)  $ be the quadratic exponential sum in $ \mathbb{P} $  defined by
$$ t\in\mathbb{P}\longmapsto \sum_{deg Y= k-1}\sum_{deg Z=  s-1}E(tYZ) \in \mathbb{Z}. $$
\small
  Then  \[ g(t)
 = \left\{\begin{array}{ccc}
   2^{s+k-j-2} &  if  & r(D_{(s-1) \times (k-1) }(t)) = r(D_{s \times ( k-1)}(t)) = r(D_{(s-1) \times k }(t)) = r(D_{s \times k}(t)) = j, \\
  -  2^{s+k-j-2} &  if  & r(D_{(s-1) \times (k-1) }(t)) = r(D_{s \times ( k-1)}(t)) = r(D_{(s-1) \times k }(t)) =  j \;  and \;r(D_{s \times k}(t)) = j+1, \\
   0    &   if  &  otherwise,
\end{array}\right.\] \\
\normalsize
and
\begin{equation*}
 \int_{\mathbb{P}}g^{2q}(t) dt =
      2^{(s+k-2)(2q-1)}\cdot  \sum_{j=0}^{s-1} {}^{\#}\Big(\begin{array}{c | c}
           j & j \\
           \hline
           j  &  j
           \end{array} \Big)_{\mathbb{P}/\mathbb{P}_{k +s -1}}\cdot2^{-2qj}.\\
            \end{equation*}
 \end{thm}

  \begin{thm}[See Definition \ref{defn 1.2}] 
 \label{thm 3.6}
 Let $ g_{m,k}(t,\eta ) = g(t,\eta ) $ be the exponential sum in $ \mathbb{P}\times\mathbb{P} $ defined by
$$ (t,\eta ) \in  \mathbb{P}\times \mathbb{P}\longmapsto 
  \sum_{deg Y\leq k-1}\sum_{deg Z\leq m}E(tYZ)\sum_{deg U =0}E(\eta YU) \in \mathbb{Z}.$$
 Then 
 \begin{equation*}
 g(t,\eta ) =  \begin{cases}
 2^{k+m+1-  r(D_{(1+m)\times k}(t)) }  & \text{if }
   r(D_{(1+m)\times k}(t)) = r(D^{\big[\stackrel{1}{1+m}\big] \times k }(t,\eta ) ), \\
     0  & \text{otherwise },
    \end{cases}
\end{equation*}
and \\
$$\int_{\mathbb{P}}\int_{\mathbb{P}}g^q(t,\eta ) dt d\eta  =
2^{q(k+m +1) -2k - m}\sum_{i = 0}^{inf(k,1+m)}
 \sigma _{i,i}^{\big[\stackrel{1}{1+m}\big]\times k}2^{-iq}. $$
 \end{thm}
  \begin{thm}[See  Definition \ref{defn 1.2}]
\label{thm 3.7}
 Let $ f_{m,k}(t,\eta ) = f(t,\eta ) $ be the exponential sum in $ \mathbb{P}\times\mathbb{P} $ defined by
$$ (t,\eta ) \in  \mathbb{P}\times \mathbb{P}\longmapsto 
  \sum_{deg Y\leq k-1}\sum_{deg Z\leq m}E(tYZ)\sum_{deg U \leq 0}E(\eta YU) \in \mathbb{Z}.$$
Then \\
$$ f(t,\eta ) = 2^{k+m +2 -r( r(D^{\big[\stackrel{1}{1+m}\big] \times k }(t,\eta )) }$$\\
and\\
   $$\int_{\mathbb{P}}\int_{\mathbb{P}}f^q(t,\eta ) dt d\eta  =
2^{q(k+m +2) -2k - m}\sum_{i = 0}^{inf(k,2+m)}
  \Gamma _{i}^{\Big[\substack{1 \\ 1+m }\Big] \times k}2^{-iq}. $$
 \end{thm}

    \begin{thm}[See  Definition \ref{defn 1.2}]
     \label{thm 3.8}
  We have the following formula  for all   $0 \leq i\leq \inf(k,2+m) $\\
  
  \begin{align*}
    \Gamma _{i}^{\Big[\substack{1 \\ 1+m }\Big] \times k}& = ( 2^{k} -  2^{i-1})\cdot\Gamma _{i-1}^{(1+m)\times k} 
 +  2^{i}\Gamma _{i}^{(1+m)\times k}. \\
  & \\
  \end{align*}

\underline {The case  k = 2}
\begin{equation*}
 \Gamma _{i}^{\Big[\substack{1 \\ 1+m }\Big] \times 2}= \begin{cases}
1 & \text{if  } i = 0, \\
 9  &  \text{if  }    i=1, \\
2^{4 +m} - 10  & \text{if   } i = 2.
\end{cases}
\end{equation*}
\underline {The case  m = 0 , $ k\geq 2 $}\\
\begin{equation*}
 \Gamma _{i}^{\Big[\substack{1 \\ 1}\Big] \times k}= \begin{cases}
 1  & \text{if  } i = 0, \\
 3\cdot(2^k-1)   &  \text{if  }    i=1, \\
2^{2k}- 3\cdot2^k + 2      & \text{if   } i = 2.
\end{cases}
\end{equation*}
\underline {The case  m = 1 , $ k\geq 3 $}\\
\begin{equation*}
 \Gamma _{i}^{\Big[\substack{1 \\ 1+ 1}\Big] \times k}= \begin{cases}
 1  & \text{if  } i = 0, \\
2^k + 5  &  \text{if  }    i=1, \\
11\cdot(2^k - 1)   & \text{if   } i = 2, \\
2^{2k+1}- 3\cdot2^{k+2} + 2^4    & \text{if   } i = 3.
\end{cases}
\end{equation*}
\underline {The case $3\leq  k \leq 1+m $}\\
\begin{equation*}
 \Gamma _{i}^{\Big[\substack{1 \\ 1+m }\Big] \times k}= \begin{cases}
1 & \text{if  } i = 0, \\
2^k +5  &  \text{if  }    i=1, \\
3\cdot2^{k+ 2i -4 } + 21\cdot2^{3i -5} & \text{if   } 2\leq i\leq k-1, \\
2^{2k +m} - 5\cdot2^{3k -5} & \text{if   } i = k.
\end{cases}
\end{equation*}
 \underline {The case $ 2\leq m\leq k-2 $}
\begin{equation*}
 \Gamma _{i}^{\Big[\substack{1 \\ 1+m }\Big] \times k}= \begin{cases}
1 & \text{if  } i = 0, \\
2^k +5  &  \text{if  }    i=1, \\
3\cdot2^{k+ 2i -4 } + 21\cdot2^{3i -5} & \text{if   } 2\leq i\leq  m, \\
11\cdot[2^{k+2m-2} - 2^{3m -2}]    & \text{if   } i = m +1,\\
2^{2k+m}  - 3\cdot2^{k+2m} +2^{3m+1}  & \text{if   } i = m +2.
\end{cases}
\end{equation*}
\end{thm}
 \begin{thm}[See Definitions \ref{defn 1.1}, \ref{defn 1.2}] 
  \label{thm 3.9}
Let $ \Gamma _{i}^{\Big[\substack{n \\ 1+m }\Big] \times k}$ denote the number of matrices of the form
 $\left[{A\over B}\right] $of rank i such that  A is a 
$(1+m)\times k$ persymmetric matrix and B is  a $ n\times k $ matrix  over $ \mathbb{F}_{2},$ and where $\Gamma _{i}^{(1+m)\times k}$ denotes
   the number of $ (1+m)\times k $ persymmetric matrices over $ \mathbb{F}_{2}$ of rank i.\\
     Then  $\Gamma _{i}^{\Big[\substack{n \\ 1+m }\Big] \times k}$ expressed as a linear combination of the $\Gamma _{i-j}^{(1+m)\times k}$ is equal to
   \begin{equation*}
 \sum_{j= 0}^{n}2^{(n-j)\cdot(i-j)} a_{j}^{(n)}\prod_{l=1}^{j}(2^{k}- 2^{i-l})\cdot
 \Gamma _{i-j}^{(1+m)\times k} \quad for \quad 0\leq i\leq inf(k,n+m +1)
 \end{equation*}
 where 
 \begin{equation*}
 a_{j}^{(n)} = \sum_{s =0}^{j-1} (-1)^{s}\prod_{l=0}^{j-(s+1)}{2^{n+1}- 2^{l}\over 2^{j-s}-2^{l}}\cdot2^{s(n-j) +{s(s+1)\over 2}}
+ (-1)^{j}\cdot 2^{jn - {j(j-1)\over 2}}  \quad for \quad  1\leq j\leq n-1.
   \end{equation*}
   We set  
   \begin{align*}
   a_{0}^{(n)} & =   a_{n}^{(n)} = 1    \\ and\quad
   \Gamma _{i-j}^{(1+m)\times k} & = 0 \quad if \quad  i-j \notin \{0,1,2,\ldots, inf(k,1+m)\}.
   \end{align*}
  \end{thm}

\begin{cor}
\label{cor 3.10}We have the following formulas for n =  1,2,3,4,5 :\\
$ \Gamma _{i}^{\Big[\substack{1 \\ 1+m }\Big] \times k}=2^{i}\Gamma _{i}^{(1+m)\times k}+
 (2^{k}-2^{i-1})\cdot\Gamma _{i-1}^{(1+m)\times k}\quad for\quad 0\leq i\leq inf(k,2+m), $\vspace{0.5 cm}\\
  $ \Gamma _{i}^{\Big[\substack{2 \\ 1+m }\Big] \times k}=2^{2i}\Gamma _{i}^{(1+m)\times k}+
     3\cdot2^{i-1}(2^{k}-2^{i-1})\cdot\Gamma _{i-1}^{(1+m)\times k} \\
     +(2^{k}-2^{i-1})(2^{k}-2^{i-2})\cdot
      \Gamma _{i-2}^{(1+m)\times k}
       \quad for\quad 0\leq i\leq inf(k,3+m), $\vspace{0.5 cm}\\
        $ \Gamma _{i}^{\Big[\substack{3 \\ 1+m }\Big] \times k}=2^{3i}\Gamma _{i}^{(1+m)\times k}+
     7\cdot2^{(i-1)2}(2^{k}-2^{i-1})\cdot\Gamma _{i-1}^{(1+m)\times k}\\
     +   7\cdot2^{i-2}(2^{k}-2^{i-1})(2^{k}-2^{i-2})\cdot \Gamma _{i-2}^{(1+m)\times k} \\
      +(2^{k}-2^{i-1})(2^{k}-2^{i-2})(2^{k}-2^{i-3}) \Gamma _{i-3}^{(1+m)\times k}
       \quad for\quad 0\leq i\leq inf(k,4+m), $\vspace{0.5 cm}\\
        $ \Gamma _{i}^{\Big[\substack{4 \\ 1+m }\Big] \times k}=2^{4i}\Gamma _{i}^{(1+m)\times k}+
     15\cdot2^{(i-1)3}(2^{k}-2^{i-1})\cdot\Gamma _{i-1}^{(1+m)\times k}\\
     +   35\cdot2^{2i-4}(2^{k}-2^{i-1})(2^{k}-2^{i-2})\cdot\Gamma _{i-2}^{(1+m)\times k}\\
      + 15\cdot2^{i-3}(2^{k}-2^{i-1})(2^{k}-2^{i-2})(2^{k}-2^{i-3}) \Gamma _{i-3}^{(1+m)\times k}\\
       + (2^{k}-2^{i-1})(2^{k}-2^{i-2})(2^{k}-2^{i-3})(2^{k}-2^{i-4}) \Gamma _{i-4}^{(1+m)\times k}\\
        for\quad 0\leq i\leq inf(k,5+m), $\vspace{0.5 cm}\\
       $ \Gamma _{i}^{\Big[\substack{5 \\ 1+m }\Big] \times k}=2^{5i}\Gamma _{i}^{(1+m)\times k}+
     31\cdot2^{(i-1)4}(2^{k}-2^{i-1})\cdot\Gamma _{i-1}^{(1+m)\times k}\\
     +   155\cdot2^{3i-6}(2^{k}-2^{i-1})(2^{k}-2^{i-2})\cdot\Gamma _{i-2}^{(1+m)\times k}\\
      + 155\cdot2^{2i-6}(2^{k}-2^{i-1})(2^{k}-2^{i-2})(2^{k}-2^{i-3}) \Gamma _{i-3}^{(1+m)\times k}\\
       +31\cdot2^{i-4}(2^{k}-2^{i-1})(2^{k}-2^{i-2})(2^{k}-2^{i-3})(2^{k}-2^{i-4}) \Gamma _{i-4}^{(1+m)\times k} \\
       + (2^{k}-2^{i-1})(2^{k}-2^{i-2})(2^{k}-2^{i-3})(2^{k}-2^{i-4})(2^{k}-2^{i-5}) \Gamma _{i-5}^{(1+m)\times k}\\
        for\quad 0\leq i\leq inf(k,6+m). $\vspace{0.5 cm}\\
    \end{cor}  
   
\begin{thm}
\label{thm 3.11}(See Section \ref{sec 1})
 Let $ \displaystyle  f_{m,k}(t,\eta_{1},\eta _{2},\ldots,\eta _{n} ) $  be the exponential sum  in $ \mathbb{P}^{n+1} $ defined by\\
  $ \displaystyle (t,\eta_{1},\eta _{2},\ldots,\eta _{n} )\in \mathbb{P}^{n+1}\longrightarrow \\
  \sum_{deg Y\leq k-1}\sum_{deg Z\leq m}E(tYZ)\sum_{deg U_{1}\leq  0}E(\eta_{1} YU_{1})
  \sum_{deg U_{2} \leq 0}E(\eta _{2} YU_{2}) \ldots \sum_{deg U_{n} \leq 0} E(\eta _{n} YU_{n}). $\vspace{0.5 cm}\\
  Set $$(t,\eta_{1},\eta _{2},\ldots,\eta _{n} ) =
  \big(\sum_{i\geq 1}\alpha _{i}T^{-i}, \sum_{i\geq 1}\beta  _{1i}T^{-i},\ldots, \sum_{i\geq 1}\beta  _{ni}T^{-i}) \in\mathbb{P}^{n+1}.   $$           
   Then
  $$ f_{m,k}(t,\eta_{1},\eta _{2},\ldots,\eta _{n} ) = 
  2^{k+m +n+1-r( D^{\Big[\stackrel{n}{1+m}\Big] \times k }(t,\eta_{1},\eta _{2},\ldots,\eta _{n} ))} $$
 where 
$$   D^{\Big[\stackrel{n}{1+m}\Big] \times k }(t,\eta_{1},\eta _{2},\ldots,\eta _{n} ) $$
 denotes  the following  $(1+n+m)\times k $ matrix
   $$  \left ( \begin{array} {cccccc}
\alpha _{1} & \alpha _{2} & \alpha _{3} &  \ldots & \alpha _{k-1}  &  \alpha _{k} \\
\alpha _{2 } & \alpha _{3} & \alpha _{4}&  \ldots  &  \alpha _{k} &  \alpha _{k+1} \\
\vdots & \vdots & \vdots   &  \ldots  & \vdots  &  \vdots \\
\vdots & \vdots & \vdots    &  \ldots & \vdots  &  \vdots \\
\alpha _{1+m} & \alpha _{2+m} & \alpha _{3+m} & \ldots  &  \alpha _{k+m-1} &  \alpha _{k+m}  \\
\hline
\beta  _{11} & \beta  _{12} & \beta  _{13} & \ldots  &  \beta_{1 k-1} &  \beta _{ 1k}  \\
\beta  _{21} & \beta  _{22} & \beta  _{23} & \ldots  &  \beta_{2 k-1} &  \beta _{2k}\\
\vdots & \vdots & \vdots   &  \ldots  & \vdots  &  \vdots \\
\vdots & \vdots & \vdots    &  \ldots & \vdots  &  \vdots \\
\beta  _{n1} & \beta  _{n2} & \beta  _{n3} & \ldots  &  \beta_{n k-1} &  \beta _{nk}
 \end{array}  \right). $$\vspace{0.5 cm}
Then the number denoted by $ R_{q}(n,k,m) $ of solutions \\
 $(Y_1,Z_1,U_{1}^{(1)},U_{2}^{(1)}, \ldots,U_{n}^{(1)}, Y_2,Z_2,U_{1}^{(2)},U_{2}^{(2)}, 
\ldots,U_{n}^{(2)},\ldots  Y_q,Z_q,U_{1}^{(q)},U_{2}^{(q)}, \ldots,U_{n}^{(q)}   ) $ \vspace{0.5 cm}\\
 of the polynomial equations  \vspace{0.5 cm}
  \[\left\{\begin{array}{c}
 Y_{1}Z_{1} +Y_{2}Z_{2}+ \ldots + Y_{q}Z_{q} = 0  \\
   Y_{1}U_{1}^{(1)} + Y_{2}U_{1}^{(2)} + \ldots  + Y_{q}U_{1}^{(q)} = 0  \\
    Y_{1}U_{2}^{(1)} + Y_{2}U_{2}^{(2)} + \ldots  + Y_{q}U_{2}^{(q)} = 0\\
    \vdots \\
   Y_{1}U_{n}^{(1)} + Y_{2}U_{n}^{(2)} + \ldots  + Y_{q}U_{n}^{(q)} = 0 
 \end{array}\right.\]
    satisfying the degree conditions \\
                   $$  degY_i \leq k-1 , \quad degZ_i \leq m ,
                   \quad degU_{j}^{i} \leq 0 , \quad  for \quad 1\leq j\leq n  \quad 1\leq i \leq q $$ \\
  is equal to the following integral over the unit interval in $ \mathbb{K}^{n+1} $
 $$ \int_{\mathbb{P}^{n+1}} f_{m,k}^{q}(t,\eta_{1},\eta _{2},\ldots,\eta _{n} )dt d\eta_{1}d\eta _{2}\ldots d\eta _{n}. $$
  Observing that $ f_{m,k}(t,\eta_{1},\eta _{2},\ldots,\eta _{n} )$ is constant on cosets of $ \mathbb{P}_{k+m}\times\mathbb{P}_{k}^{n}, $\;
  the above integral is equal to 
$$  2^{q(k+m +n+1) -(n+1)k - m}\sum_{i = 0}^{inf(k,n+1+m)}
   \Gamma _{i}^{\left[\stackrel{n}{1+m}\right]\times k}2^{-iq} = R_{q}(n,k,m)$$
\end{thm}
\begin{example}
\label{example 3.12}
The number $R_{q}(0,k,m)$  of solutions
 $(Y_1,Z_1, \ldots,Y_q,Z_q) $  of the polynomial equation
                        $$ Y_1Z_1 +  Y_2Z_2 + \ldots + Y_qZ_q = 0  $$ \\
  satisfying the degree conditions \\
                   $$  degY_i \leq k-1 , \quad degZ_i \leq m \leq k-1 \quad for  1\leq i \leq q. $$ \\
 is equal to the following integral   \\ 
  \begin{align*}
& \int_{\mathbb{P}}\Big[ \sum_{deg Y\leq k-1}\sum_{deg Z\leq m}E(tYZ) \Big]^{q}dt
  = 2^{(q-1)(k+m+1) +1}\sum_{i=0}^{1+m}\Gamma_{i}^{(1+m)\times k}2^{-qi}  \\
 &  = \left\{\begin{array}{ccc}
   2^k + 2^{1+m} -1 & if & q=1,\\
   2^{2k} + 3\cdot(m+1)\cdot 2^{k+m} & if & q=2, \\
   2^{(q-1)(k+m+1) +1}\left[1+ 3\frac{1-2^{(2-q)m}}{2^{q}- 2^{2}}+
 (2^{k+m}-2^{2m})2^{-q(1+m)}\right]
 & if & 3\leq q.
  \end{array}\right. 
  \end{align*}
\end{example}

\begin{example}
\label{example 3.13}
The number  $ \Gamma _{i}^{\Big[\substack{1 \\ 1+2 }\Big] \times 3}$ of rank i matrices of the form \\
  $$  \left ( \begin{array} {ccc}
\alpha _{1} & \alpha _{2} & \alpha _{3}  \\
\alpha _{2 } & \alpha _{3} & \alpha _{4} \\
\alpha _{3} & \alpha _{4} & \alpha _{5} \\
\hline
\beta  _{1} & \beta  _{2} & \beta  _{3} \\
 \end{array}  \right) $$
is equal to 

\begin{equation*}
  \begin{cases}
1 & \text{if  } i = 0, \\
13 &  \text{if  }    i=1, \\
66  & \text{if   } i = 2, \\
176  & \text{if   } i = 3.
\end{cases}
\end{equation*}
\vspace{0.1 cm}

 The number $R_{q}(1,3,2)$ of solutions 
 $(Y_1,Z_1,U_{1}, \ldots,Y_q,Z_q,U_{q}) $  of the polynomial equations
   \[\left\{\begin{array}{c}
 Y_{1}Z_{1} +Y_{2}Z_{2}+ \ldots + Y_{q}Z_{q} = 0  \\
   Y_{1}U_{1} + Y_{2}U_{2} + \ldots  + Y_{q}U_{q} = 0
 \end{array}\right.\]
  satisfying the degree conditions \\
                   $$  degY_i \leq 2 , \quad degZ_i \leq 2 ,\quad degU_{i}\leq 0 \quad for \quad 1\leq i \leq q $$ 
  is equal to the following integral \\
   \begin{align*}
  & \int_{\mathbb{P}}\int_{\mathbb{P}}\big[  \sum_{deg Y\leq 2}\sum_{deg Z\leq 2}E(tYZ)\sum_{deg U \leq 0}E(\eta YU) \big]^{q}dtd\eta  
  =  2^{7q-8}\sum_{i = 0}^{3}\Gamma _{i}^{\Big[\substack{1 \\ 1+2 }\Big] \times 3}2^{-iq}\\
  & \\
  & = 2^{4q-8}\cdot\big[2^{3q} +13\cdot2^{2q} +66\cdot2^{q} +176 \big].
   \end{align*}
\end{example}

\begin{example}
\label{example 3.14}
 The number  $ \Gamma _{i}^{\Big[\substack{5 \\ 1+2 }\Big] \times 4}$ of rank i matrices of the form \\
  $$  \left ( \begin{array} {cccc}
\alpha _{1} & \alpha _{2} & \alpha _{3} & \alpha _{4} \\
\alpha _{2 } & \alpha _{3} & \alpha _{4} & \alpha _{5}\\
\alpha _{3} & \alpha _{4} & \alpha _{5} & \alpha _{6} \\
\hline
\beta  _{11} & \beta  _{12} & \beta  _{13} & \beta_{14}  \\
\beta  _{21} & \beta  _{22} & \beta  _{23} & \beta_{24}  \\
\beta  _{31} & \beta  _{32} & \beta  _{33} & \beta_{34}  \\
\beta  _{41} & \beta  _{42} & \beta  _{43} & \beta_{44}  \\
\beta  _{51} & \beta  _{52} & \beta  _{53} & \beta_{54} 
 \end{array}  \right) $$
is equal to 

\begin{equation*}
  \begin{cases}
1 & \text{if  } i = 0, \\
561 &  \text{if  }    i=1, \\
65670 & \text{if   } i = 2, \\
3731208 & \text{if   } i = 3, \\
63311424   & \text{if   } i = 4. 
\end{cases}
\end{equation*}
\vspace{0.1 cm}

 The number $R_{3}(5,4,2)$ of solutions \\
 
  $(Y_1,Z_1,U_{1}^{(1)},U_{2}^{(1)},U_{3}^{(1)},U_{4}^{(1)},U_{5}^{(1)}, Y_2,Z_2,U_{1}^{(2)},U_{2}^{(2)}, 
U_{3}^{(2)},U_{4}^{(2)},U_{5}^{(2)}, Y_3,Z_3,U_{1}^{(3)},U_{2}^{(3)},U_{3}^{(3)},U_{4}^{(3)},U_{5}^{(3)}   ) $ \\

 of the polynomial equations
  \[\left\{\begin{array}{c}
 Y_{1}Z_{1} +Y_{2}Z_{2} + Y_{3}Z_{3} = 0,  \\
   Y_{1}U_{1}^{(1)} + Y_{2}U_{1}^{(2)}  + Y_{3}U_{1}^{(3)} = 0,  \\
 Y_{1}U_{2}^{(1)} + Y_{2}U_{2}^{(2)}  + Y_{3}U_{2}^{(3)} = 0, \\
 Y_{1}U_{3}^{(1)} + Y_{2}U_{3}^{(2)}  + Y_{3}U_{3}^{(3)} = 0, \\
 Y_{1}U_{4}^{(1)} + Y_{2}U_{4}^{(2)}  + Y_{3}U_{4}^{(3)} = 0, \\
  Y_{1}U_{5}^{(1)} + Y_{2}U_{5}^{(2)}  + Y_{3}U_{5}^{(3)} = 0, 
 \end{array}\right.\]
   satisfying the degree conditions \\
                   $$  degY_i \leq 3 , \quad degZ_i \leq 2 ,\quad degU_{i}\leq 0 \quad  for \quad 1\leq j\leq 5  \quad 1\leq i \leq 3 $$ \\ 
    is equal  to the following integral over the unit interval in $ \mathbb{K}^{6} $\\
   \begin{align*}
  & \int_{\mathbb{P}^{6}} f_{2,4}^{3}(t,\eta_{1},\eta _{2},\eta _{3},\eta _{4},\eta _{5} )dt d\eta_{1}d\eta _{2}d\eta _{3}d\eta _{4}d\eta _{5}
  =  2^{10}\cdot \sum_{i = 0}^{4}\Gamma _{i}^{\Big[\substack{5 \\ 1+2 }\Big] \times 4}2^{-i3}= 24413824.
   \end{align*}

 \end{example}

  \section{Exponential sums formulas  on $ \mathbb{K} $}
   \label{sec 4}
 In this section we compute exponential quadratic sums in  $\mathbb{P}$ and show that they only depend on rank properties of some 
 associated persymmetric matrices.  
 The following propositions are proved in [2]. \\
\begin{prop}
\label{prop 4.1}
The following holds
\begin{itemize}
\item  For every rational integer j , the measure of  $\mathbb{P}_{j} \; is \; 2^{-j}.$
\item   For every \; $ A \in \mathbb{F}_{2}[T]  $,  \;E(A) = 1 .
\item   For $ u \in \mathbb{K} \quad \nu(u)\geq  2 \Rightarrow  $ \; E(u) = 1.
\end{itemize}
\end{prop}
\begin{prop}
\label{prop 4.2}
Let j be a rational integer and \; $ u \in \mathbb{K} $,\; then 
  \[ \int_{\mathbb{P}_{j}}E(ut)dt   =  \left\{ \begin{array}{ccc}
                         2^{-j}    &  if  &  \nu(u) > - j,    \\
                          0  &       &  otherwise.
                          \end{array}\right.\]
\end{prop}
 \begin{prop}
\label{prop 4.3}
 Let j be a rational integer and let \; $ u \in \mathbb{K},$ \;then
    \[ \sum_{deg B\leq  j}E(Bu)  =  \left\{ \begin{array}{ccc}
                         2^{j+1}    &  if  &   \nu\left(\left\{u\right\}\right) > j +1,  \\
                          0  &       &  otherwise.
                          \end{array}\right.\]
\end{prop}
\begin{lem}
\label{lem 4.4}
Let $Y \in\mathbb{F}_{2}[T], $\;  then\\
  \[ \sum_{deg Z\leq s-1} E(tYZ)   =  \left\{ \begin{array}{ccc}
                         2^{s}    &  if  &   \nu\left(\left\{tY\right\}\right) > s,\\
                          0  &       &  otherwise.
                          \end{array}\right.\]
\end{lem}
 \begin{proof}
 
 Lemma \ref{lem 4.4} follows from  Proposition \ref{prop 4.3} with u = tY and j = s -1.
\end{proof}
\begin{lem}
\label{lem 4.5}
Let $ t \in\mathbb{P}\; and \; Y \in \mathbb{F}_{2}[T]  ,\; deg Y\leq k-1, $\; then\\
$$\nu\left(\left\{tY\right\}\right) > s
 \Longleftrightarrow  Y\in \ker D_{s\times k}(t). $$
\end{lem}
 \begin{proof}
 Let  $ Y=\sum_{j=0}^{k-1}\gamma _{j}T^j     \; , \gamma _{j}\in \mathbb{F}_{2}, $\; then\\
 $$ tY =(\sum_{i\geq 1}\alpha _{i}T^{-i})(\sum_{j=0}^{k-1} \gamma  _{j}  T^{j})
  = \sum_{i\geq 1}\sum_{j=0}^{k-1}\alpha_{i}\gamma_{j}T^{-(i-j)} $$ and  \\
 $$\left\{tY\right\} = \left(\sum_{i=1}^{k}
 \alpha_{i}\gamma_{i-1}\right)T^{-1}
+\left(\sum_{i=2}^{k+1}\alpha_{i}\gamma_{i-2}\right)T^{-2} +
\ldots + \left(\sum_{i= s}^{k+s-1}\alpha_{i}\gamma_{i-s}\right)T^{- s} + \ldots $$
Therefore $  \mathcal{\nu}(\left\{tY\right\}) \geq s  $  if and only if 
  \def\BA{ \begin{array}{ccccccc}}
  \def\BB{ \begin{array}{c}}
  \def\EA{\end{array}}
  \begin{displaymath}
  \left[\BA
\alpha  _{1} & \alpha_{2} &\alpha _{3} & \alpha_{4} &\alpha _{5}  & \ldots  & \alpha_{k}\\
\alpha _{2} &\alpha _{3} &\alpha _{4} & \alpha_{5} & \alpha_{6}  & \ldots  & \alpha_{k+1} \\
 \vdots & \vdots & \vdots &\vdots & \vdots & \ldots                   & \vdots \\
 \vdots & \vdots & \vdots & \vdots & \vdots  & \ldots                   &\vdots \\
 \alpha _{s} &\alpha _{s+1} & \alpha_{s+2} &\alpha _{s+3} & \alpha_{s+4} &\ldots  &\alpha _{k+s -1} \\  \EA\right ]
 \left[ \BB \gamma _{0}\\ \gamma _{1}\\ \gamma _{2}\\ \gamma _{3}\\ \gamma _{4}
  \\ \vdots  \\\gamma _{k-1} \EA\right] =
  \left[ \BB 0 \\ 0 \\ 0 \\  0 \\ 0  \\ \vdots \\  0
    \EA\right]
 \end{displaymath} \\
\begin{displaymath} \Longleftrightarrow
D_{s\times k}(t) \left(\begin{array}{c}
\gamma _{0}\\
\gamma _{1}\\
\vdots \\
\gamma _{k-1}\\
\end{array}\right) =\left(\begin{array}{c}
0\\
0\\
\vdots\\
0\\
\end{array}\right)\Longleftrightarrow
Y= \left(\begin{array}{c}
\gamma _{0}\\
\gamma _{1}\\
\vdots \\
\gamma _{k-1}\\
\end{array}\right)   \in Ker \, D_{s\times k}(t).
\end{displaymath}
\end{proof}

\begin{lem}
\label{lem 4.6}
Let $ t \in\mathbb{P} $ and let  q be  a rational integer $\geq 1, $ then \\
$$   h_{s,k}(t) = h(t) = \sum_{deg Y\leq k-1}\sum_{deg Z\leq s-1}E(tYZ) =
2^{k + s - r\big(D_{s \times k}(t)\big)}$$ and 
$$ \int_{\mathbb{P}}h^q(t) = 2^{(q-1)\cdot(k+s)+1}\sum_{i=0}^{s}\Gamma _{i}^{s\times k}2^{-qi}.$$
\end{lem}
 \begin{proof}
 By Lemma \ref{lem 4.4} and Lemma \ref{lem 4.5} we obtain \\
$$ \sum_{deg Y\leq k-1}\sum_{deg Z\leq s-1}E(tYZ) =
 2^{s}\sum_{deg Y\leq k-1 \atop{Y \in Ker D_{s \times k}(t)}}1 = 
 2^{s}\cdot2^{k -r\big( D_{s \times k}(t)\big)}=
 2^{s +k -r\big( D_{s \times k}(t)\big)}.$$
 Observing that h(t) is constant on cosets of $ \mathbb{P}_{k+s-1}, $ we obtain by integrating $ h^q(t) $
 on the unit interval of $ \mathbb{K}$\vspace{0.5 cm}
 \begin{align*}
 \int_{\mathbb{P}}h^{q}(t) dt & =
\sum_{t \in \mathbb{P}/\mathbb{P}_{k+s-1}}2^{[s+k -r(D_{s\times k}(t))]\cdot q}
\int_{\mathbb{P}_{k+s-1}}dt  \\ & =
 \sum_{i = 0}^{s} \Gamma _{i}^{s\times k}\cdot 2^{(s+k-i)\cdot q}\int_{\mathbb{P}_{k+s-1}}dt \\
  &
=\sum_{i = 0}^{s} \Gamma _{i}^{s\times k}\cdot 2^{(s+k-i)\cdot q}\cdot 2^{-(k+s-1)} \\
& =
 2^{(q-1)(k+s) +1}\sum_{i=0}^{s}\Gamma_{i}2^{-qi}. 
 \end{align*}
  \end{proof}
\begin{lem}
\label{lem 4.7}
Let  $ t \in\mathbb{P}$ and $g_{s,k}(t) =  g(t)= \sum_{deg Y=k-1}\sum_{deg Z=s-1}E(tYZ), $ \;then \\
 $$g^{2}(t) = g_{1}(t)\cdot g_{2}(t) \quad   where  $$ \\
 $g_{1}(t) =  \sum_{deg Y\leq k-2}\sum_{deg Z= s-1}E(tYZ) \; and \quad  g_{2}(t) =  \sum_{deg Y= k -1}\sum_{deg Z \leq s-2}E(tYZ). $
\end{lem}
 \begin{proof} Observe that 
 $$ g^2(t) =
   \sum_{deg Y= k-1}\sum_{deg Y_{1}=k-1}\sum_{deg Z= s-1}
   \sum_{deg Z_{1} = s-1}E(tYZ) E(tY_{1}Z_{1}). $$ \\
 Set      \[
  \left\{\begin{array}{cc}
Y +Y_{1} = Y_{2}       &  deg Y_{2} \leq k-2,  \\
 Y_{1} = Y_{3}     &   deg Y_{3} = k-1, \\
 Z +Z_{1} = Z_{2}       &  deg Z_{2} \leq s-2,  \\
 Z_{1} = Z_{3}     &   deg Z_{3} = s-1.
\end{array}\right.\] \\
Then we obtain \\
\begin{eqnarray*}
 g^2(t) & =  & \sum_{degY_{2}\leq k-2} \sum_{degY_{3}= k-1}
  \sum_{degZ_{2}\leq s-2} \sum_{degZ_{3}= s -1}E(t(Y_{2}+ Y_{3}) Z_{3})
  E(tY_{3}(Z_{2} + Z_{3})) \\
  & = & \sum_{degY_{2}\leq k-2} \sum_{degY_{3}= k-1}
  \sum_{degZ_{2}\leq s-2} \sum_{degZ_{3}= s -1}E(tY_{2}Z_{3})E(tY_{3}Z_{2}) \\
&   =  &  \big[ \sum_{deg Y\leq k-2}\sum_{deg Z=s-1}E(tYZ) \big]\cdot
    \big[ \sum_{deg Y= k-1}\sum_{deg Z\leq s-2}E(tYZ) \big] \\
    & = &  g_{1}(t)\cdot g_{2}(t).
  \end{eqnarray*}
  \end{proof}
  \begin{lem}
\label{lem 4.8}
 Let $ t \in\mathbb{P}, $ \; then
  \begin{equation*}
g_{1}(t)  =   \begin{cases}
 2^{k+s-2 - r(D_{(s-1)\times (k-1)}(t))}  & \text{if }
   r(D_{(s-1)\times( k-1)}(t)) = r(D_{s\times( k-1)}(t)), \\
     0  & \text{otherwise },
    \end{cases}
\end{equation*} and \\
 \begin{equation*}
g_{2}(t)  =   \begin{cases}
 2^{k+s-2 - r(D_{(s-1)\times (k-1)}(t))}  & \text{if      }
   r(D_{(s-1)\times( k-1)}(t)) = r(D_{(s-1)\times k}(t)), \\
     0  & \text{otherwise    }.
    \end{cases}
\end{equation*}
 \end{lem}
 \begin{proof} Observe that 
 \begin{eqnarray*}
 g_{1}(t) & = &  \sum_{deg Y\leq k-2}\sum_{deg Z= s-1}E(tYZ) \\
 & = & \sum_{deg Y\leq k-2}\sum_{deg Z\leq s-1}E(tYZ)  -
  \sum_{deg Y\leq k-2}\sum_{deg Z \leq s-2}E(tYZ)\\ \vspace{0.5 cm}
  & = & h_{s,k-1}(t) - h_{s-1,k-1}(t).
\end{eqnarray*} \\
By Lemma \ref{lem 4.6} we obtain \\
$$ g_{1}(t) = 2^{k + s -1 - r\big(D_{s \times(k-1)}(t)\big)} -
  2^{k + s -2 - r\big(D_{(s-1) \times(k-1) }(t)\big)}.$$
 In the same way we obtain \\
$$  g_{2}(t) = 2^{k + s -1 - r\big(D_{(s-1) \times k}(t)\big)} -
  2^{k + s -2 - r\big(D_{(s-1) \times(k-1) }(t)\big)}$$\\
   and Lemma \ref{lem 4.8} is proved.
 \end{proof}
  \begin{lem}We have the following equivalences
\label{lem 4.9}
\begin{eqnarray*}
g(t) \neq 0 & \Longleftrightarrow  & g^{2}(t) \neq  0 \\
& \Longleftrightarrow
&  g_{1}(t)\cdot g_{2}(t) \neq 0 \\
& \Longleftrightarrow &
  r\big(D_{(s-1) \times(k-1) }(t)\big) = r\big(D_{s \times(k-1) }(t)\big) =
   r\big(D_{(s-1) \times k }(t)\big)
\end{eqnarray*}
\end{lem}
 \begin{proof}
 Lemma \ref{lem 4.9} is obtained by combining Lemma \ref{lem 4.7} and Lemma \ref{lem 4.8}
\end{proof}
  \begin{lem}
\label{lem 4.10}
Let $ t\in\mathbb{P},$\; then $ g_{s,k}(t) = g(t) $ is given by \\
\begin{itemize}
\item  $ 2^{s+k-j-2} \quad if \quad r(D_{(s-1) \times (k-1) }(t)) = r(D_{s \times ( k-1)}(t)) = r(D_{(s-1) \times k }(t)) = r(D_{s \times k}(t)) = j $
\item  $ -  2^{s+k-j-2} \quad if \quad  r(D_{(s-1) \times (k-1) }(t)) = r(D_{s \times ( k-1)}(t)) = r(D_{(s-1) \times k }(t)) =  j \;  and \;r(D_{s \times k}(t)) = j+1 $
\item     0    \quad    otherwise.
\end{itemize} 
\end{lem}
 \begin{proof} We see that 
  \begin{equation*}
 g(t)  =   \sum_{deg Y= k-1}\sum_{deg Z= s-1}E(tYZ)  =  \sum_{deg Y\leq k-1}\sum_{deg Z= s-1}E(tYZ)  -
  \sum_{deg Y\leq k-2}\sum_{deg Z = s-1}E(tYZ).
\end{equation*} \\
By Lemma \ref{lem 4.9}\\
$$g(t) \neq 0  \Longleftrightarrow
  r\big(D_{(s-1) \times(k-1) }(t)\big) = r\big(D_{s \times(k-1) }(t)\big) =
   r\big(D_{(s-1) \times k }(t)\big). $$
   We consider the following  two cases: \\
   First case\\
   $$  r\big(D_{(s-1) \times(k-1) }(t)\big) = r\big(D_{s \times(k-1) }(t)\big) =
   r\big(D_{(s-1) \times k }(t)\big) =   r\big(D_{s \times k }(t)\big) = j. $$
   Second case
  $$  r\big(D_{(s-1) \times(k-1) }(t)\big) = r\big(D_{s \times(k-1) }(t)\big) =
   r\big(D_{(s-1) \times k }(t)\big) = j $$ \\
  and    $$   r\big(D_{s \times k }(t)\big) = j +1. $$
  By the results of lemma \ref{lem 4.8} we obtain:  \\
 In the first case
$$ \sum_{deg Y\leq k-1}\sum_{deg Z= s-1}E(tYZ)  -
  \sum_{deg Y\leq k-2}\sum_{deg Z = s-1}E(tYZ) $$\\
  $$  = 2^{k +s -1 -j} - 2^{k +s -2 -j}=  2^{k +s -2 -j}.$$\\
 In the second case
 $$ \sum_{deg Y\leq k-1}\sum_{deg Z= s-1}E(tYZ)  -
  \sum_{deg Y\leq k-2}\sum_{deg Z = s-1}E(tYZ) $$\\
 $$  = 0  - 2^{k +s -2 -j}=  - 2^{k +s -2 -j}$$
 and Lemma \ref{lem 4.10} is proved.
\end{proof}
\section{Rank properties of a partition of persymmetric matrices.} [See Definition \ref{defn 1.1} ].\quad
 \label{sec 5}
  In this section we establish rank properties of a partition of the following $ s\times k $ persymmetric matrix,
   obtained by drawing in the matrix  a horizontal line between the (s -1)th  and the sth row and drawing a vertical line between 
  the (k-1)th and the kth column  
 \begin{equation*}
 t = \sum_{i\geq 1}\alpha _{i}T^{-1}\in \mathbb{P} \longrightarrow 
 \left ( \begin{array} {ccccc|c}
\alpha _{1} & \alpha _{2} & \alpha _{3} &  \ldots & \alpha _{k-1}  &  \alpha _{k} \\
\alpha _{2 } & \alpha _{3} & \alpha _{4}&  \ldots  &  \alpha _{k} &  \alpha _{k+1} \\
\vdots & \vdots & \vdots   &  \vdots  & \vdots  & \vdots \\
\vdots & \vdots & \vdots    &  \vdots & \vdots   & \vdots   \\
\alpha _{s-1} & \alpha _{s} & \alpha _{s+1} & \ldots  &  \alpha _{k+s-3} &  \alpha _{k+s-2}  \\
\hline
\alpha _{s} & \alpha _{s+1} & \alpha _{s+2} & \ldots  &  \alpha _{k+s-2} &  \alpha _{k+s-1}  \\
\end{array}  \right)  \in Persym_{s\times k}(\mathbb{F}_{2})
\end{equation*}

     \begin{lem}
\label{lem 5.1}
Let $ q\in\mathbb{N},$\;  then\\
$$\int_{\mathbb{P}}g^{2q+1}(t)dt = 0. $$
 \end{lem}
 \begin{proof}
 Recall that by Lemma \ref{lem 4.7}
$$ g_{s,k}(t) =  g(t) = \sum_{deg Y= k-1}\sum_{deg Z= s-1}E(tYZ), $$
then $$ \int_{\mathbb{P}}g^{2q+1}(t)dt $$
$$ = \sum_{degY_{1} = k-1}\sum_{deg Z_{1}= s-1}\ldots
\sum_{degY_{2q +1} = k-1}\sum_{deg Z_{2q +1}= s-1}
\int_{\mathbb{P}}\prod_{i=1}^{2q+1}E(tY_{i}Z_{i})dt $$
$$ = \sum_{degY_{1} = k-1}\sum_{deg Z_{1}= s-1}\ldots
\sum_{degY_{2q +1} = k-1}\sum_{deg Z_{2q +1}= s-1}
\int_{\mathbb{P}}E\Big(t\sum_{i =1}^{2q +1}Y_{i}Z_{i}\Big)dt. $$\\
 By Proposition \ref{prop 4.2} with j = 0
 $$\int_{\mathbb{P}}
 E(t(Y_{1}Z_{1}+ Y_{2}Z_{2} +\ldots  +Y_{2q +1}Z_{2q +1}))dt  = 0  \Longleftrightarrow \sum_{i =1}^{2q +1}Y_{i}Z_{i} = 0 $$
 but 2q +1 is odd, therefore  $ deg\Big(\sum_{i =1}^{2q +1}Y_{i}Z_{i}\Big) = s +k -2. $
  \end{proof}
  
  \begin{lem}
\label{lem 5.2}
 Let $ j \in\mathbb{N} \; such\; that \; 0\leq j\leq s-1, $\; then \\
  \begin{align*}
{}^{\#}\Big(\begin{array}{c | c}
           j & j \\
           \hline
           j  &  j
           \end{array} \Big)_{\mathbb{P}/\mathbb{P}_{k +s -1}}  =
           {}^{\#}\Big(\begin{array}{c | c}
           j & j \\
           \hline
           j  &  j +1
           \end{array} \Big)_{\mathbb{P}/\mathbb{P}_{k +s -1}}
            \end{align*} (see Definition \ref{defn 1.1}).
            \end{lem}
  \begin{proof}
   We define
$$\mathbb{A} = \{ t \in \mathbb{P}
\mid r(D_{(s-1)\times( k-1)}(t)) =  r(D_{(s-1) \times k}(t))
= r(D_{s\times( k-1)}(t)) =  r(D_{s \times k}(t))  \}$$
and \\
$$\mathbb{B} = \{ t \in \mathbb{P}
\mid r(D_{(s-1)\times( k-1)}(t)) =  r(D_{(s-1) \times k}(t))  = r(D_{s\times( k-1)}(t)) ,
  r(D_{s \times k}(t)) =  r(D_{(s-1)\times( k-1)}(t))+1\}.$$
 By lemma \ref{lem 4.10} we have by observing that g(t) is constant on cosets of $ \mathbb{P}_{k+s-1} $
$$ \int_{\mathbb{P}}g^{2q+1}(t)dt
 = \int_{  \mathbb{A}} g^{2q+1}(t)dt
    +  \int_{ \mathbb{B}} g^{2q+1}(t)dt $$
$$ =\sum_{t \in \mathbb{A}\bigcap\big(\mathbb{P}/\mathbb{P}_{k+s-1}\big) }
 2^{(s+k-2-  r(D_{(s-1)\times( k-1)}(t)))(2q-1)}\int_{\mathbb{P}_{s+k-1}}dt $$
$$  +  \sum_{t \in \mathbb{B}\bigcap\big(\mathbb{P}/\mathbb{P}_{k+s-1}\big)}
- 2^{(s+k-2-  r(D_{(s-1)\times( k-1)}(t)))(2q+1)}\int_{\mathbb{P}_{s+k-1}}dt $$
  \begin{align*}
=   \sum_{j=0}^{s-1}  2^{(s+k-2-j)(2q+1)}\cdot {}^{\#}\Big(\begin{array}{c | c}
           j & j \\
           \hline
           j  &  j
           \end{array} \Big)_{\mathbb{P}/\mathbb{P}_{k +s -1}} \int_{\mathbb{P}_{k+s-1}}dt \\
             -    \sum_{j=0}^{s-1}  2^{(s+k-2-j)(2q+1)}\cdot {}^{\#}\Big(\begin{array}{c | c}
           j & j \\
           \hline
           j  &  j+1
           \end{array} \Big)_{\mathbb{P}/\mathbb{P}_{k +s -1}} \int_{\mathbb{P}_{k+s-1}}dt
              \end{align*}
By lemma \ref{lem 5.1},\; for all  $ q \in\mathbb{N} $
 \begin{align*}
  \sum_{j=0}^{s-1}  2^{(-j)(2q+1)}\cdot\Big( {}^{\#}\Big(\begin{array}{c | c}
           j & j \\
           \hline
           j  &  j
           \end{array} \Big)_{\mathbb{P}/\mathbb{P}_{k +s -1}}
               -   {}^{\#}\Big(\begin{array}{c | c}
           j & j \\
           \hline
           j  &  j+1
           \end{array} \Big)_{\mathbb{P}/\mathbb{P}_{k +s -1}}\Big) = 0. 
              \end{align*}
This completes the proof.
\end{proof}
  \begin{lem} [See Definition \ref{defn 1.1}]
\label{lem 5.3}
 Let $ j \in\mathbb{N} \; such\; that \; 0\leq j\leq s-2, $\; then\\
  \begin{equation}
  \label{eq 5.1}
{}^{\#}\Big(\begin{array}{c | c}
           j & j+1 \\
           \hline
           j +1 &  j +1
           \end{array} \Big)_{\mathbb{P}/\mathbb{P}_{k +s -1}}  = 0, 
            \end{equation}
 \begin{equation}
 \label{eq 5.2}
 \left.{\begin{array}{cc}
    r(D_{(s-1)\times (k-1)}(t)) =  r(D_{s \times (k-1)}(t)) = j     \\
 r(D_{(s-1) \times (k-2)}^{2}(t)) = j-1  
\end{array}}\right\} \Longrightarrow  r(D_{s \times (k-2)}^{2}(t)) = j-1.  
\end{equation}
   \end{lem}
 \begin{proof}
 We consider the following partition of the matrix $ D_{s\times k}(t) $
  $$ \left ( \begin{array} {c|cccc|c}
\alpha _{1} & \alpha _{2} & \alpha _{3} &  \ldots & \alpha _{k-1}  &  \alpha _{k} \\
\alpha _{2 } & \alpha _{3} & \alpha _{4}&  \ldots  &  \alpha _{k} &  \alpha _{k+1} \\
\vdots & \vdots & \vdots   &  \vdots  & \vdots  & \vdots \\
\vdots & \vdots & \vdots    &  \vdots & \vdots & \vdots \\
\alpha _{s-1} & \alpha _{s} & \alpha _{s+1} & \ldots  &  \alpha _{k+s-3} &  \alpha _{k+s-2}  \\
\hline
\alpha _{s} & \alpha _{s+1} & \alpha _{s+2} & \ldots  &  \alpha _{k+s-2} &  \alpha _{k+s-1}  \\
\end{array}  \right) $$ \\
Obviously we have $ r(D_{(s-1) \times k}(t)) =  r(D_{s \times k}(t)) \Longrightarrow
r(D_{(s-1)\times( k-1)}(t)) =  r(D_{s\times( k-1)}(t)) $ and so \eqref{eq 5.1} is proved.
We remark that this result is true for any matrix.
To prove \eqref{eq 5.2} we consider the matrix obtained by the matrix $ D_{s\times k}(t) $ 
by deleting the first column and replacing the last column by the first one :
  $$ \left ( \begin{array} {cccc|c}
 \alpha _{2} & \alpha _{3} &  \ldots & \alpha _{k-1}  &  \alpha _{1} \\
 \alpha _{3} & \alpha _{4}&  \ldots  &  \alpha _{k} &  \alpha _{2} \\
\vdots & \vdots & \vdots   &  \vdots  & \vdots    \\
\vdots & \vdots & \vdots    &  \vdots & \vdots   \\
 \alpha _{s} & \alpha _{s+1} & \ldots  &  \alpha _{k+s-3} &  \alpha _{s-1}  \\
\hline
 \alpha _{s+1} & \alpha _{s+2} & \ldots  &  \alpha _{k+s-2} &  \alpha _{s}  \\
\end{array}  \right). $$ \\ 
Now \eqref{eq 5.2} follows from  \eqref{eq 5.1} with  j replaced by  j-1.
\end{proof}
 \begin{lem}
\label{lem 5.4}
 For every  $ j \in\mathbb{N} \; such\; that \; 0\leq j\leq s-2 $ we have \\
  \begin{align*}
{}^{\#}\Big(\begin{array}{c | c}
           j & j +1 \\
           \hline
           j  &  j + 1
           \end{array} \Big)_{\mathbb{P}/\mathbb{P}_{k +s -1}}  =
           {}^{\#}\Big(\begin{array}{c | c}
           j & j \\
           \hline
           j +1 &  j +1
           \end{array} \Big)_{\mathbb{P}/\mathbb{P}_{k +s -1}} = 0
            \end{align*} (see Definition \ref{defn 1.1}).
            \end{lem}
 \begin{proof}
Proof by contradiction.\\
Assume on the contrary  that there exists $ j_{0} \in [0,s-2] $  such that
 \begin{equation}
\label{eq 5.3}
{}^{\#}\Big(\begin{array}{c | c}
           j_{0} & j_{0} +1 \\
           \hline
           j_{0}  &  j_{0} + 1
           \end{array} \Big)_{\mathbb{P}/\mathbb{P}_{k +s -1}}  > 0.
 \end{equation}
 We are going to show that\\
 \small
 \begin{equation*}
{}^{\#}\Big(\begin{array}{c | c}
           j_{0} & j_{0} +1 \\
           \hline
           j_{0}  &  j_{0} + 1
           \end{array} \Big)_{\mathbb{P}/\mathbb{P}_{k +s -1}}  > 0  \Longrightarrow
{}^{\#}\Big(\begin{array}{c | c}
           j_{0}-1 & j_{0}  \\
           \hline
           j_{0}-1  &  j_{0}
           \end{array} \Big)_{\mathbb{P}_{1}/\mathbb{P}_{k +s -1}}  > 0  \Longrightarrow
 \ldots  \Longrightarrow
{}^{\#}\Big(\begin{array}{c | c}
           0 & 1  \\
           \hline
           0 &  1
           \end{array} \Big)_{\mathbb{P}_{j}/\mathbb{P}_{k +s -1}}  > 0,
 \end{equation*}
 \normalsize
  which obviously contradicts
 \begin{equation*}
 {}^{\#}\Big(\begin{array}{c | c}
           0 & 1  \\
           \hline
           0 &  1
           \end{array} \Big)_{\mathbb{P}_{j}/\mathbb{P}_{k +s -1}}  = 0.
 \end{equation*}
  By   \eqref{eq 5.3} there exists $ t_{0}\in \mathbb{P}/\mathbb{P}_{k +s -1} $ such that \\
$$ r(D_{(s-1)\times( k-1)}(t_{0})) =  r(D_{s\times( k-1)}(t_{0})) = j_{0} $$ and 
$$  r(D_{(s-1) \times k}(t_{0})) =   r(D_{s \times k}(t_{0})) = j_{0}+1. $$
 Consider the following partition of  the matrix $ D_{s\times k}(t_{0})$\\
  $$  \left ( \begin{array} {c|cccc|c}
\alpha _{1} & \alpha _{2} & \alpha _{3} &  \ldots & \alpha _{k-1}  &  \alpha _{k} \\
\hline
\alpha _{2 } & \alpha _{3} & \alpha _{4}&  \ldots  &  \alpha _{k} &  \alpha _{k+1} \\
\vdots & \vdots & \vdots   &  \vdots  & \vdots  &\vdots \\
\vdots & \vdots & \vdots    &  \vdots & \vdots  & \vdots \\
\alpha _{s-1} & \alpha _{s} & \alpha _{s+1} & \ldots  &  \alpha _{k+s-3} &  \alpha _{k+s-2}  \\
\hline
\alpha _{s} & \alpha _{s+1} & \alpha _{s+2} & \ldots  &  \alpha _{k+s-2} &  \alpha _{k+s-1}  \\
\end{array}  \right). $$ \\
 By \eqref{eq 5.3},  
  $ r(D_{(s-1)\times( k-1)}^{2}(t_{0}))\leq  j_{0} $ and $ r(D_{(s-1) \times k}(t_{0})) = j_{0}+1 $, it follows that  $ r(D_{(s-1)\times( k-1)}^{2}(t_{0}))= j_{0}$.\\
  
    Let $ a_{1},a_{2},\ldots,a_{k}$ denote the  columns of \quad $ D_{(s-1)\times k}(t_{0}), $  that is
\begin{displaymath}
a_{i} = \left(\begin{array}{c}
\alpha_{i}\\
\alpha_{i+1}\\
\vdots \\
\alpha_{i+s-2}\\
\end{array}\right)
\end{displaymath} 
  Since  $ r(D_{(s-1) \times (k-1)}^2(t_{0})) = j_{0}\; and \; r(D_{(s-1)\times k}(t_{0})) = j_{0} +1 $, we have 
$  a_{1}\notin span\left\{a_{2},a_{3},\ldots,a_{k}\right\} $,
therefore  $ r(D_{(s-1) \times( k-2)}^2(t_{0})) = j_{0}-1$.\vspace{0.5 cm}\\
Now  $ r(D_{(s-1) \times( k-2)}^2(t_{0})) = j_{0}-1,
\quad  r(D_{(s-1) \times( k-1)}(t_{0})) = j_{0} \quad and \quad r(D_{s \times( k-1)}(t_{0})) = j_{0}$ .\\
From \eqref{eq 5.2} we get $ r(D_{s \times( k-2)}^2(t_{0})) =  j_{0}-1. $
\vspace{0.5 cm} Thus we obtain \\
 $r(D_{(s-1) \times( k-2)}^2(t_{0})) = r(D_{s \times( k-2)}^2(t_{0})) = j_{0}-1 \quad and \quad
r(D_{(s-1) \times( k-1)}^2(t_{0})) =  r(D_{s \times( k-1)}^2(t_{0})) = j_{0}$.\vspace{0.5 cm}\\
 Consider now the matrix 
$$ D_{s \times(k-1)}^2(t_{0})    = \left ( \begin{array} {ccccc |c}
\alpha _{2} & \alpha _{3} & \alpha _{4} &  \ldots & \alpha _{k-1}  & \alpha _{k} \\
\alpha _{3 } & \alpha _{4} & \alpha _{5} &  \ldots & \alpha _{k}  & \alpha _{k+1} \\
\vdots & \vdots   &  \vdots  & \ldots  & \vdots & \vdots \\
\vdots & \vdots & \vdots  &   \ldots  & \vdots & \vdots \\
\alpha _{s } & \alpha _{s+1} & \alpha _{s+2} & \ldots & \alpha _{k+s-3} & \alpha _{k+s-2}  \\
\hline \\
\alpha _{s +1} & \alpha _{s+2} & \alpha _{s+2} & \ldots & \alpha _{k+s-2} & \alpha _{k+s-1}  \\
\end{array}  \right ).$$
We  get  \begin{align*}
{}^{\#}\Big (\begin{array}{c| c}
       j_{0} -1& j_{0} \\
       \hline
       j_{0}-1 & j_{0}
\end{array}\Big)_{ \mathbb{P}_{1}/\mathbb{P}_{k+s-1}} > 0.
\end{align*} \\
We now have  proved that \\
 \begin{equation*}
{}^{\#}\Big(\begin{array}{c | c}
           j_{0} & j_{0} +1 \\
           \hline
           j_{0}  &  j_{0} + 1
           \end{array} \Big)_{\mathbb{P}/\mathbb{P}_{k +s -1}}  > 0  \Longrightarrow
{}^{\#}\Big(\begin{array}{c | c}
           j_{0}-1 & j_{0}  \\
           \hline
           j_{0}-1  &  j_{0}
           \end{array} \Big)_{\mathbb{P}_{1}/\mathbb{P}_{k +s -1}}  > 0.
\end{equation*}
 We repeat this procedure and obtain after finitely many steps \\
 \begin{align*}
{}^{\#}\Big (\begin{array}{c| c}
        0 & 1 \\
       \hline
       0 & 1
\end{array}\Big)_{ \mathbb{P}_{j_{0}}/\mathbb{P}_{k+s-1}} > 0.
\end{align*} \\
From    \begin{displaymath}
 D_{s \times( k-j_{0})}^{j_{0}+1}(t_{0})    = \left ( \begin{array} {ccccc |c}
\alpha _{j_{0}+1} & \alpha _{j_{0}+2} & \alpha _{j_{0}+3} &  \ldots & \alpha _{k-1}  & \alpha _{k} \\
\alpha _{j_{0}+2} & \alpha _{j_{0}+3} & \alpha _{j_{0}+4} &  \ldots & \alpha _{k}  & \alpha _{k+1} \\
\vdots & \vdots   &  \vdots  & \ldots  & \vdots & \vdots \\
\vdots & \vdots & \vdots  &   \ldots  & \vdots & \vdots \\
\alpha _{s+j_{0}-1 } & \alpha _{s+j_{0}} & \alpha _{s+j_{0}+1} & \ldots & \alpha _{k+s-3} & \alpha _{k+s-2}  \\
\hline\\
\alpha _{s +j_{0}} & \alpha _{s+j_{0}+1} & \alpha _{s+j_{0}+1} & \ldots & \alpha _{k+s-2} & \alpha _{k+s-1}  \\
\end{array}  \right ),
\end{displaymath} \\  we  obviously get 
\begin{align*}
{}^{\#}\Big (\begin{array}{c| c}
        0 & 1 \\
       \hline
       0 & 1
\end{array}\Big)_{ \mathbb{P}_{j_{0}}/\mathbb{P}_{k+s-1}} = 0.
\end{align*} \\
 Arguing in the same way with the transpose of $D_{s\times k}(t) ,$ observing that the rank of a matrix is equal to the rank of its transpose,
  yields  $  {}^{\#}\Big(\begin{array}{c | c}
           j & j \\
           \hline
           j +1 &  j +1
           \end{array} \Big)_{\mathbb{P}/\mathbb{P}_{k +s -1}} = 0 $.                                                                                    

\end{proof}
\begin{lem}
\label{lem 5.5}
We have \\
\begin{equation*}
   \Gamma_{i}^{s \times (k-1)} =   \Gamma_{i}^{(s-1) \times k}\quad  \text{if}
   \quad  0\leq i \leq s-2, \quad s\leq k. 
\end{equation*}
 \end{lem}
 \begin{proof}
 Consider  the Matrix $ D_{s\times k}(t) $
$$\left(
   \begin{array} {ccccc | c}
\alpha _{1} & \alpha _{2} & \alpha _{3}&  \ldots  & \alpha _{k-1}  & \alpha _{k}\\
\alpha _{2 } & \alpha _{3} & \alpha _{4}&  \ldots  & \alpha _{k} & \alpha _{k+1} \\
\vdots & \vdots & \vdots   &  \vdots  &  \vdots   &  \vdots \\
\vdots & \vdots & \vdots    &  \vdots & \vdots   &  \vdots \\
\alpha _{s-1 } & \alpha _{s} & \alpha _{s+1} & \ldots  & \alpha _{k+s-3}  & \alpha _{k+s-2}  \\
\hline \\
\alpha _{s } & \alpha _{s+1} & \alpha _{s+2} & \ldots & \alpha _{k+s-2}  & \alpha _{k+s-1} \\
\end{array}
\right). $$
 By Lemmas \ref{lem 5.2}, \ref{lem 5.3}  and \ref{lem 5.4} we have \\
\Small 
\begin{align*}
2\cdot\Gamma_{i}^{s \times (k-1)} &  = &
  {}^{\#}\Big (\begin{array}{c| c}
       i  & i \\
       \hline
       i & i
\end{array}\Big)_{ \mathbb{P}/\mathbb{P}_{k+s-1}} + 
 {}^{\#}\Big (\begin{array}{c| c}
       i  & i \\
       \hline
       i & i+1
\end{array}\Big)_{ \mathbb{P}/\mathbb{P}_{k+s-1}} +
 {}^{\#}\Big (\begin{array}{c| c}
       i  & i +1 \\
       \hline
       i & i +1
\end{array}\Big)_{ \mathbb{P}/\mathbb{P}_{k+s-1}}  +
 {}^{\#}\Big (\begin{array}{c| c}
       i-1  & i -1 \\
       \hline
       i & i - 1
\end{array}\Big)_{ \mathbb{P}/\mathbb{P}_{k+s-1}} \\
& + & {}^{\#}\Big (\begin{array}{c| c}
       i-1  & i -1 \\
       \hline
       i & i 
\end{array}\Big)_{ \mathbb{P}/\mathbb{P}_{k+s-1}} +
 {}^{\#}\Big (\begin{array}{c| c}
       i-1  & i  \\
       \hline
       i & i 
\end{array}\Big)_{ \mathbb{P}/\mathbb{P}_{k+s-1}}+
 {}^{\#}\Big (\begin{array}{c| c}
       i-1  & i  \\
       \hline
       i & i +1
\end{array}\Big)_{ \mathbb{P}/\mathbb{P}_{k+s-1}}\\
& = & 2\cdot{}^{\#}\Big (\begin{array}{c| c}
       i  & i \\
       \hline
       i & i
\end{array}\Big)_{ \mathbb{P}/\mathbb{P}_{k+s-1}} 
+{}^{\#}\Big (\begin{array}{c| c}
        i-1 & i \\
       \hline
       i & i+1
\end{array}\Big)_{ \mathbb{P}/\mathbb{P}_{k+s-1}} \\
\intertext{and}
2\cdot\Gamma_{i}^{(s-1) \times k} & = &
 {}^{\#}\Big (\begin{array}{c| c}
        i & i \\
       \hline
       i & i
\end{array}\Big)_{ \mathbb{P}/\mathbb{P}_{k+s-1}}
+{}^{\#}\Big (\begin{array}{c| c}
        i & i \\
       \hline
       i & i+1
\end{array}\Big)_{ \mathbb{P}/\mathbb{P}_{k+s-1}}
+{}^{\#}\Big (\begin{array}{c| c}
        i & i  \\
       \hline
       i+1 & i+1
\end{array}\Big)_{ \mathbb{P}/\mathbb{P}_{k+s-1}}\\
& + & {}^{\#}\Big (\begin{array}{c| c}
        i-1 & i \\
       \hline
       i-1 & i
\end{array}\Big)_{ \mathbb{P}/\mathbb{P}_{k+s-1}}
 + {}^{\#}\Big (\begin{array}{c| c}
        i-1 & i \\
       \hline
       i & i
\end{array}\Big)_{ \mathbb{P}/\mathbb{P}_{k+s-1}}
+{}^{\#}\Big (\begin{array}{c| c}
        i-1 & i \\
       \hline
       i & i+1
\end{array}\Big)_{ \mathbb{P}/\mathbb{P}_{k+s-1}}\\
& = & 2\cdot {}^{\#}\Big (\begin{array}{c| c}
        i & i \\
       \hline
       i & i
\end{array}\Big)_{ \mathbb{P}/\mathbb{P}_{k+s-1}}
+{}^{\#}\Big (\begin{array}{c| c}
        i-1 & i \\
       \hline
       i & i+1
\end{array}\Big)_{ \mathbb{P}/\mathbb{P}_{k+s-1}}.
\end{align*}
\end{proof}
\normalsize
\begin{lem}
\label{lem 5.6}
We have (see Definition \ref{defn 1.1})
\begin{equation*}
 \int_{\mathbb{P}}g^{2q}(t) dt =
      2^{(s+k-2)(2q-1)}\cdot  \sum_{j=0}^{s-1} {}^{\#}\Big(\begin{array}{c | c}
           j & j \\
           \hline
           j  &  j
           \end{array} \Big)_{\mathbb{P}/\mathbb{P}_{k +s -1}}\cdot2^{-2qj}.
            \end{equation*}
\end{lem}
\begin{proof}
  From Lemma \ref{lem 5.2},  by observing that g(t) is constants on cosets of  $\mathbb{P}_{k+s -1}, $ \; we obtain
 \begin{align*}
  \int_{\mathbb{P}}g^{2q}(t) dt =
      \sum_{j=0}^{s-1}  2^{(s+k-2-j)2q}\cdot {}^{\#}\Big(\begin{array}{c | c}
           j & j \\
           \hline
           j  &  j
           \end{array} \Big)_{\mathbb{P}/\mathbb{P}_{k +s -1}} \int_{\mathbb{P}_{k+s-1}}dt  & \\
            +   \sum_{j=0}^{s-1}  2^{(s+k-2-j)2q}\cdot {}^{\#}\Big(\begin{array}{c | c}
           j & j \\
           \hline
           j  &  j+1
           \end{array} \Big)_{\mathbb{P}/\mathbb{P}_{k +s -1}} \int_{\mathbb{P}_{k+s-1}}dt & = \\
            2\cdot \sum_{j=0}^{s-1}  2^{(s+k-2-j)2q}\cdot {}^{\#}\Big(\begin{array}{c | c}
           j & j \\
           \hline
           j  &  j
           \end{array} \Big)_{\mathbb{P}/\mathbb{P}_{k +s -1}}2^{-(k+s-1)} & = \\
            2^{(s+k-2)(2q-1)}\cdot  \sum_{j=0}^{s-1} {}^{\#}\Big(\begin{array}{c | c}
           j & j \\
           \hline
           j  &  j
           \end{array} \Big)_{\mathbb{P}/\mathbb{P}_{k +s -1}}\cdot2^{-2qj}
             \end{align*} which prove Lemma \ref{lem 5.6}.                                                                             
  \end{proof}
  \section{Exponential  sums in  $\mathbb{K}\times\mathbb{K} $}
 \label{sec 6}[See  Definition \ref{defn 1.2}].
 In this section we compute exponential sums in  $\mathbb{P}^{2}, $ associated to a matrix of the form  $\big[{A\over b_{-}}\big],$ 
 where A is a $ (m+1)\times k $ persymmetic matrix over  $ \mathbb{F}_{2} $ and $ b_{-} $ a $1\times k $ matrix with entries in  $ \mathbb{F}_{2}.$
 \begin{lem}
\label{lem 6.1}
Set   $  \eta  = \sum_{j\geq 1}\beta _{j}T^{-j} \in \mathbb{P}\;
  and \;Y = \sum_{i=1}^{k-1}\delta _{i}T^i\in \mathbb{F}_{2}[T] ,\quad degY\leq k-1.  $ \\
  Then 
  \begin{equation*}
E(\eta Y) = \begin{cases}
1 & \text{if        }                      \sum_{j=1}^{k-1}\beta _{j}\delta _{j} = 0, \\
-1 & \text{if        }                       \sum_{j=1}^{k-1}\beta _{j}\delta _{j} = 1.
\end{cases}
\end{equation*}
  \end{lem}
\begin{proof}
By lemmas \ref{lem 4.4},\;\ref{lem 4.5} with $ t\longrightarrow \eta $ and s = 1 we have \\
$$\nu\left(\left\{\eta Y\right\}\right) > 1
 \Longleftrightarrow  Y\in \ker D_{1\times k}(\eta )
 \Longleftrightarrow   \sum_{j=1}^{k-1}\beta _{j}\delta _{j} = 0,  $$
$$\nu\left(\left\{\eta Y\right\}\right) = 1
 \Longleftrightarrow   \sum_{j=1}^{k-1}\beta _{j}\delta _{j} = 1.  $$
\end{proof}
\begin{lem}
\label{lem 6.2}
Let $ (t,\eta ) \in  \mathbb{P}\times \mathbb{P} $ and set 
$$      g_{m,k}(t,\eta ) =  g(t,\eta ) = \sum_{deg Y\leq k-1}\sum_{deg Z\leq m}E(tYZ)\sum_{deg U =0}E(\eta YU). $$  \\
Then \\
  \begin{equation*}
 g(t,\eta ) = \begin{cases}
 2^{k+m+1-  r(D_{(1+m)\times k}(t)) }  & \text{if }
   r(D_{(1+m)\times k}(t)) = r(D^{\big[\stackrel{1}{1+m}\big] \times k }(t,\eta ) ), \\
     0  & \text{otherwise}.
    \end{cases}
\end{equation*}
\end{lem}
\begin{proof}
 Consider the matrix [see Definition \ref{defn 1.2}] \\
 $$  D^{\Big[\stackrel{1}{1+m}\Big] \times k }(t,\eta )  = \left ( \begin{array} {cccccc}
\alpha _{1} & \alpha _{2} & \alpha _{3} &  \ldots & \alpha _{k-1}  &  \alpha _{k} \\
\alpha _{2 } & \alpha _{3} & \alpha _{4}&  \ldots  &  \alpha _{k} &  \alpha _{k+1} \\
\vdots & \vdots & \vdots   &  \vdots  & \vdots  &  \vdots \\
\vdots & \vdots & \vdots    &  \vdots & \vdots  &  \vdots \\
\alpha _{1+m} & \alpha _{2+m} & \alpha _{3+m} & \ldots  &  \alpha _{k+m-1} &  \alpha _{k+m}  \\
\hline
\beta  _{1} & \beta  _{2} & \beta  _{3} & \ldots  &  \beta_{k-1} &  \beta _{k}  \\
\end{array}  \right). $$ \\
 By Lemmas \ref{lem 4.4}, \ref{lem 4.5} and \ref{lem 6.1} we obtain \\
\begin{align*}
 g(t,\eta ) &  = 2^{1+m}\sum_{deg Y\leq k-1\atop {Y\in ker D_{(1+m)\times k}(t) }}E(\eta Y) \\
& = 2^{1+m}\cdot\Big[\sum_{deg Y\leq k-1\atop {Y\in ker D^{\big[\stackrel{1}{1+m}\big] \times k }(t,\eta )}} 1
 -   \sum_{\substack{deg Y\leq k-1 \\
 Y\in ker D_{(1+m)\times k}(t) \\
 \sum_{j = 0}^{k -1}\beta _{j}\delta _{j}= 1}}1 \Big] \\
 & =  2^{1+m}\cdot\Big[\sum_{deg Y\leq k-1\atop {Y\in ker D^{\big[\stackrel{1}{1+m}\big] \times k }(t,\eta )}} 1  -
  \big(\sum_{\substack{deg Y\leq k-1 \\
 Y\in ker D_{(1+m)\times k}(t)}}1  -   \sum_{\substack{deg Y\leq k-1 \\
Y\in ker D^{\big[\stackrel{1}{1+m}\big] \times k }(t,\eta )}}1 \big )\Big] \\
& =   2^{1+m}\cdot\Big[2\cdot\sum_{deg Y\leq k-1\atop {Y\in ker D^{\big[\stackrel{1}{1+m}\big] \times k }(t,\eta )}} 1  -  \sum_{\substack{deg Y\leq k-1 \\
 Y\in ker D_{(1+m)\times k}(t)}}1 \Big]  \\
 & = 2^{1+m}\cdot\Big[2\cdot 2^{k - r( D^{\big[\stackrel{1}{1+m}\big] \times k }(t,\eta ))}- 2^{k - r(D_{(1+m)\times k}(t)) }\Big].
\end{align*}
\end{proof}
\begin{lem}
\label{lem 6.3}
Set  $  h(t)  = \sum_{deg \leq k-1}\sum_{deg Z \leq m}E(tYZ) = h_{m+1,k}(t) $ and let $ q\geq 2 $ be an integer, then \\
            $$ g^q(t,\eta ) = g(t,\eta )\cdot h^{q-1}(t). $$
\end{lem}
\begin{proof}Obviously we get 
$$ g^2(t,\eta ) =\big[ 2^{1+m}\sum_{deg Y_{1}\leq k-1\atop {Y_{1}\in ker D_{(1+m)\times k}(t) }}E(\eta Y_{1}) \big]\cdot
 \big[ 2^{1+m}\sum_{deg Y_{2}\leq k-1\atop {Y_{2}\in ker D_{(1+m)\times k}(t) }}E(\eta Y_{2}) \big]. $$
Define
\[
\left\{\begin{array}{cc}
Y_{1} + Y_{2} = Y_{3}  &  deg Y_{3}\leq k-1 \\
              Y_{1} = Y_{4}   &   deg Y_{4}\leq k-1.
\end{array}\right.\]
Then we obtain  \\
\begin{align*}
 g^2(t,\eta ) & = 2^{2(1+m)}\sum_{deg Y_{1}\leq k-1\atop {Y_{1}\in ker D_{(1+m)\times k}(t) }}
 \sum_{deg Y_{2}\leq k-1\atop {Y_{2}\in ker D_{(1+m)\times k}(t) }}E(\eta (Y_{1}+ Y_{2}))  \\
 & = \big[ 2^{1+m}\sum_{deg Y_{4}\leq k-1\atop {Y_{4}\in ker D_{(1+m)\times k}(t) }} 1 \big]\cdot
 \big[ 2^{1+m}\sum_{deg Y_{3}\leq k-1\atop {Y_{3}\in ker D_{(1+m)\times k}(t) }} E(\eta Y_{ 3})\big] \\
 & = h(t)\cdot g(t,\eta ).
\end{align*}
By recurrence on q we get   $g^q(t,\eta ) = g(t,\eta )\cdot h^{q-1}(t). $
\end{proof}
\begin{lem}
\label{lem 6.4} We have
$$\int_{\mathbb{P}}\int_{\mathbb{P}}g^q(t,\eta ) dt d\eta  =
2^{q(k+m +1) -2k - m}\sum_{i = 0}^{inf(k,1+m)}
 \sigma _{i,i}^{\big[\stackrel{1}{1+m}\big]\times k}2^{-iq}. $$
 \end{lem}
 \begin{proof}  Recall  Definition \ref{defn 1.2},  then Lemma \ref{lem 6.2} gives,  by observing that $ g(t,\eta ) $ is
 constant on cosets of \\
  $ \mathbb{P}_{k+m}\times\mathbb{P}_{k}.$
 \begin{align*}
\int_{\mathbb{P}}\int_{\mathbb{P}}g^q(t,\eta ) dt d\eta & =
 \sum_{(t,\eta )\in \mathbb{P}/\mathbb{P}_{k+m}\times \mathbb{P}/\mathbb{P}_{k}\atop
 { r(D_{(1+m)\times k}(t)) = r(D^{\big[\stackrel{1}{1+m}\big] \times k }(t,\eta ))}}
 2^{q(k+m+1- r(D_{(1+m)\times k}(t))}\int_{\mathbb{P}_{k+m}}dt \int_{\mathbb{P}_{k}}d\eta \\
 & = \sum_{i = 0}^{inf(k,1+m)}
 \sigma _{i,i}^{\big[\stackrel{1}{1+m}\big]\times k} 2^{q(k+m+1- i)}2^{-(k+m)}2^{-k}.
  \end{align*}
 \end{proof}

\begin{lem} Equally
\label{lem 6.5}
$$ \int_{\mathbb{P}}\int_{\mathbb{P}}g^q(t,\eta ) dt d\eta  =
2^{q(k+m +1) -2k - m}\sum_{i = 0}^{inf(k,1+m)}2^{i}\Gamma _{i}^{(1+m)\times k}
2^{-iq}. $$
\end{lem}
\begin{proof}
From  Lemma  \ref{lem 6.3} and  Lemma \ref{lem 4.6} with s = m+1  we get by Fubini's theorem  \\ 
 \begin{align*}
\int_{\mathbb{P}}\int_{\mathbb{P}}g^q(t,\eta ) dt d\eta  & =
\int_{\mathbb{P}}\int_{\mathbb{P}}g(t,\eta )h^{q-1}(t) dt d\eta   =\int_{\mathbb{P}}h^{q-1}(t)\big(\int_{\mathbb{P}}g(t,\eta )d\eta \big)dt 
 =  2^{1+m}\int_{\mathbb{P}}h^{q-1}(t)dt \\
 & = 2^{1+m}\cdot 
\sum_{t \in \mathbb{P}/\mathbb{P}_{k+m}}2^{[k+m+1 -r(D_{(m+1)\times k}(t))]\cdot (q-1)}
\int_{\mathbb{P}_{k+m}}dt  \\ & = 2^{1+m}\cdot 
 \sum_{i = 0}^{inf(k,1+m)} \Gamma _{i}^{(m+1)\times k}\cdot 2^{(m+1+k-i)\cdot (q-1)}\int_{\mathbb{P}_{k+m}}dt \\
  & = 2^{1+m}\cdot2^{(q-2)(k+m+1)+1}\sum_{i=0}^{inf(k,1+m)}\Gamma _{i}^{(1+m)\times k}2^{-(q-1)i}. 
\end{align*}
\end{proof}
\begin{lem}
\label{lem 6.6}  For all $ 0\leq i\leq inf(k, 1+m) $ we have $ \sigma _{i,i}^{\big[\stackrel{1}{1+m}\big]\times k} =
 2^{i}\Gamma _{i}^{(1+m)\times k}. $ 
\end{lem}
\begin{proof}
Lemma \ref{lem 6.4}  and Lemma \ref{lem 6.5} give \\
$$\sum_{i=0}^{inf(k,1+m)}\big( \sigma _{i,i}^{\big[\stackrel{1}{1+m}\big]\times k}  -  2^{i}
\Gamma _{i}^{(1+m)\times k} )\cdot2^{-iq} = 0 \quad  for\; all\; q \geq 2. $$\\
\end{proof}
\begin{lem}
\label{lem 6.7} For all  $1\leq i\leq inf(k,2+m) $  we have\\
 $ \sigma _{i-1,i}^{\big[\stackrel{1}{1+m}\big]\times k} = (2^k - 2^{i-1})\Gamma _{i-1}^{(1+m)\times k}.$
\end{lem}
\begin{proof} Consider the matrix  \\
  $$ D^{\Big[\stackrel{1}{1+m}\Big] \times k }(t,\eta )  = \left ( \begin{array} {cccccc}
\alpha _{1} & \alpha _{2} & \alpha _{3} &  \ldots & \alpha _{k-1}  &  \alpha _{k} \\
\alpha _{2 } & \alpha _{3} & \alpha _{4}&  \ldots  &  \alpha _{k} &  \alpha _{k+1} \\
\vdots & \vdots & \vdots   &  \vdots  & \vdots  &  \vdots \\
\vdots & \vdots & \vdots    &  \vdots & \vdots  &  \vdots \\
\alpha _{1+m} & \alpha _{2+m} & \alpha _{3+m} & \ldots  &  \alpha _{k+m-1} &  \alpha _{k+m}  \\
\hline
\beta  _{1} & \beta  _{2} & \beta  _{3} & \ldots  &  \beta_{k-1} &  \beta _{k}  \\
\end{array}  \right). $$ \\
Clearly Lemma  \ref{lem 6.6} yields
\begin{align*}
2^{k}\cdot\Gamma _{i-1}^{(1+m)\times k} & = \sigma _{i-1,i-1}^{\big[\stackrel{1}{1+m}\big]\times k}+ \sigma _{i-1,i}^{\big[\stackrel{1}{1+m}\big]\times k} \\
& \Leftrightarrow  \sigma _{i-1,i}^{\big[\stackrel{1}{1+m}\big]\times k}  = 2^{k}\cdot\Gamma _{i-1}^{(1+m)\times k} - 2^{i-1}\cdot\Gamma _{i-1}^{(1+m)\times k}\\
&  \Leftrightarrow  \sigma _{i-1,i}^{\big[\stackrel{1}{1+m}\big]\times k}  = ( 2^{k} -  2^{i-1})\cdot\Gamma _{i-1}^{(1+m)\times k}    && \text{if } 1\leq i\leq inf(k,2+m)
 \end{align*}
\end{proof}
\begin{lem}
\label{lem 6.8}
For all  $0\leq i\leq inf(k,2+m) $ we have \\
 $ \Gamma _{i}^{\Big[\substack{1 \\ 1+m }\Big] \times k} = ( 2^{k} -  2^{i-1})\cdot\Gamma _{i-1}^{(1+m)\times k}  +  2^{i}\Gamma _{i}^{(1+m)\times k}. $
 \end{lem}
\begin{proof} From Lemmas \ref{lem 6.6}, \ref{lem 6.7} we get \\
 $$ \Gamma _{i}^{\Big[\substack{1 \\ 1+m }\Big] \times k}
=    \sigma _{i,i}^{\big[\stackrel{1}{1+m}\big]\times k}+ \sigma _{i-1,i}^{\big[\stackrel{1}{1+m}\big]\times k}
 = 2^{i}\Gamma _{i}^{(1+m)\times k} + ( 2^{k} -  2^{i-1})\cdot\Gamma _{i-1}^{(1+m)\times k}. $$
  \end{proof}
  \begin{lem}
\label{lem 6.9}
Set $ f(t,\eta ) =  \sum_{deg Y\leq k-1}\sum_{deg Z\leq m}E(tYZ)\sum_{deg U \leq 0}E(\eta YU) $\\
$ \quad for \quad  k\geq 1\quad and  \quad  m\geq  0 . $\\
Then \\
$$ f(t,\eta ) = 2^{k+m +2 -r( r(D^{\big[\stackrel{1}{1+m}\big] \times k }(t,\eta )) }$$\\
and \\
   $$\int_{\mathbb{P}}\int_{\mathbb{P}}f^q(t,\eta ) dt d\eta  =
2^{q(k+m +2) -2k - m}\sum_{i = 0}^{inf(k,2+m)}
  \Gamma _{i}^{\Big[\substack{1 \\ 1+m }\Big] \times k}2^{-iq}. $$
\end{lem}
\begin{proof}
By using  Lemma \ref{lem 4.4} and Lemma \ref{lem 4.5} we obtain  \\
\begin{align*}
f(t,\eta ) &  = \sum_{deg Y\leq k-1}\sum_{deg Z\leq m}E(tYZ)\sum_{deg U \leq 0}E(\eta YU) \\
& = 2^{2+m}\sum_{{ deg Y\leq k-1\atop Y\in ker D_{(1+m)\times k}(t)}\atop Y\in ker D_{1\times k}(\eta )}  1\\
& =  2^{2+m}\sum_{ deg Y\leq k-1\atop Y\in ker D^{\big[\stackrel{1}{1+m}\big] \times k }(t,\eta )} 1 =  2^{k+m +2 -r( r(D^{\big[\stackrel{1}{1+m}\big] \times k }(t,\eta )}.
 \end{align*}
 Then we have, by observing that $ f(t,\eta ) $ is constant on cosets of
  $ \mathbb{P}_{k+m}\times\mathbb{P}_{k}$
   \begin{align*}
\int_{\mathbb{P}}\int_{\mathbb{P}}f^q(t,\eta ) dt d\eta & =
 \sum_{(t,\eta )\in \mathbb{P}/\mathbb{P}_{k+m}\times \mathbb{P}/\mathbb{P}_{k}}
 2^{q(k+m+2- r( r(D^{\big[\stackrel{1}{1+m}\big] \times k }(t,\eta )) }\int_{\mathbb{P}_{k+m}}dt \int_{\mathbb{P}_{k}}d\eta \\
 & = \sum_{i = 0}^{inf(k,2+m)} \Gamma _{i}^{\Big[\substack{1 \\ 1+m }\Big] \times k} 2^{q(k+m+2- i)}2^{-(k+m)}2^{-k}.
  \end{align*}
\end{proof}

 \section{Proof of  Theorem 3.1}
\label{sec 7}
We prove Theorem \ref{thm 3.1} by  induction on s   ( recall that $ s\leq k $).\\
 Set  $ h_{s,k}(t) =  \sum_{deg Y\leq k-1}\sum_{deg Z\leq s-1}E(tYZ). $ Let s = 2 and consider the matrix $ D_{2\times k}(t), $
that is 
\begin{equation*}
 D_{2 \times k}(t)    = \left ( \begin{array} {ccccc}
\alpha _{1} & \alpha _{2} & \alpha _{3}&  \ldots  \alpha _{k} \\
\alpha _{2 } & \alpha _{3} & \alpha _{4}&  \ldots  \alpha _{k+1} \\
\end{array}  \right ).
\end{equation*}
We have
\begin{equation}
\label{eq 7.1}
\sum_{i=0}^2 \Gamma _{i}^{2 \times k} = 2^{k+1}.
\end{equation}
By Lemma \ref{lem 4.6} with s = 2 and q = 1
\begin{eqnarray}
\label{eqn 7.2}
 \int_{\mathbb{P}} h_{2,k}(t)dt  = 
\int_{\mathbb{P}} 2^{k+2-r( D_{2\times k}(t))}dt  =  2\cdot\sum_{i=0}^2 \Gamma _{i}^{2\times k}\cdot2^{-i}. 
\end{eqnarray}
 On the other hand\\
 
 \begin{equation}
 \label{eq 7.3}
  \int_{\mathbb{P}} h_{2,k}(t)dt  =  
   Card\left\{(Y,Z),degY\leq k-1,degZ\leq 1\mid Y\cdot Z = 0\right\} =  2^k +2^2 -1. 
\end{equation}
From  \eqref{eqn 7.2} and \eqref{eq 7.3} we obtain
\begin{equation}
\label{eq 7.4}
2\Gamma_{1}^{2\times k} + \Gamma_{2}^{2\times k}= 2^{k+1} +2.
\end{equation}
 From \eqref{eq 7.1} and \eqref{eq 7.4} we deduce
 \[\Gamma_{i}^{2\times k} = \left\{\begin{array}{ccc}
             1 & if & i= 0, \\
             3 & if &  i =1, \\
             2^{(k+1)} -4  & if & i = 2.
             \end{array}\right.\]
Let s = 3  and consider the matrix $ D_{3\times k}(t) $
\begin{equation*}
 D_{3 \times k}(t)    = \left ( \begin{array} {ccccc}
\alpha _{1} & \alpha _{2} & \alpha _{3}&  \ldots  \alpha _{k} \\
\alpha _{2 } & \alpha _{3} & \alpha _{4}&  \ldots  \alpha _{k+1} \\
\alpha _{3} & \alpha _{4} & \alpha _{5}&  \ldots  \alpha _{k+2} \\
\end{array}  \right ).
\end{equation*}
We have
\begin{equation}
\label{eq 7.5}
\sum_{i=0}^3 \Gamma _{i}^{3 \times k} = 2^{k+2}.
\end{equation}
By Lemma \ref{lem 5.5}\\
\begin{equation}
\label{eq 7.6}
   \Gamma_{i}^{3 \times k} =   \Gamma_{i}^{2 \times (k+1)}\;
  if \; 0\leq i \leq 1, \quad k\geq 3.
\end{equation}
Using  \eqref{eq 7.5} and  \eqref{eq 7.6}, we get  \\
\begin{equation}
\label{eq 7.7}
\Gamma _{2}^{3\times k} + \Gamma _{3}^{3\times k} = 2^{k+2} -2^2.
\end{equation}
By Lemma \ref{lem 4.6} with s = 3 and q =1
\begin{eqnarray}
\label{eqn 7.8}
 \int_{\mathbb{P}} h_{3,k}(t)dt = \int_{\mathbb{P}} 2^{k+3-r( D_{3\times k}(t))}dt 
 =  2\cdot\sum_{i=0}^3 \Gamma _{i}^{3\times k}\cdot2^{-i}. 
\end{eqnarray}
 On the other hand
 \begin{equation}
 \label{eq 7.9}
  \int_{\mathbb{P}} h_{3,k}(t)dt  =  
    Card\left\{(Y,Z),degY\leq k-1,degZ\leq 2\mid Y\cdot Z = 0\right\} =  2^k +2^3 -1. 
\end{equation}
From  \eqref{eqn 7.8} and \eqref{eq 7.9} we obtain
\begin{equation}
\label{eq 7.10}
2\Gamma_{2}^{3\times k} + \Gamma_{3}^{3\times k}= 2^{k+2} +2^3.
\end{equation}
 From  \eqref{eq 7.7} and \eqref{eq 7.10} we deduce
 \[\Gamma_{i}^{3\times k} = \left\{\begin{array}{ccc}
             1 & if & i= 0, \\
             3 & if &  i =1, \\
             12 & if & i = 2, \\
              2^{(k+2)} - 2^4  & if & i = 3.
             \end{array}\right.\]
Let s = 4. In the same way we obtain \\
 \[\Gamma_{i}^{4\times k} = \left\{\begin{array}{ccc}
             1 & if & i= 0, \\
             3 & if &  i =1, \\
             12 & if & i = 2, \\
             48 & if & i=3, \\
              2^{(k+3)} - 2^6  & if & i = 4.
             \end{array}\right.\]
 Assume that for  $ k\geq s-1$ we have
 \begin{equation*}
      (H) \quad  \Gamma_{i}^{(s-1)\times k} = \left\{\begin{array}{ccc}
             1 & if & i= 0, \\
  3\cdot2^{2(i-1)}  &   if & 1\leq i\leq s-2, \\
   2^{(k+s-2)} - 2^{2s-4}  & if & i = s-1.
             \end{array}\right.
             \end{equation*}
From (H) and Lemma \ref{lem 5.5} it follows
  \begin{equation}
   \label{eq 7.11}
    \Gamma_{i}^{s \times (k-1)} =   \Gamma_{i}^{(s-1) \times k} =  3\cdot2^{2(i-1)}\;
  if \; 0\leq i \leq s-2 \quad , k\geq s+1.
        \end{equation}
  Consider the following partition of the matrix  $  D_{s\times k}(t)$ \\
 $$  \left( 
   \begin{array} {ccccc | c}
\alpha _{1} & \alpha _{2} & \alpha _{3}&  \ldots  & \alpha _{k-1}  & \alpha _{k}\\
\alpha _{2 } & \alpha _{3} & \alpha _{4}&  \ldots  & \alpha _{k} & \alpha _{k+1} \\
\vdots & \vdots & \vdots   &  \vdots  &  \vdots   &  \vdots \\
\vdots & \vdots & \vdots    &  \vdots & \vdots   &  \vdots \\
\alpha _{s-1 } & \alpha _{s} & \alpha _{s+1} & \ldots  & \alpha _{k+s-3}  & \alpha _{k+s-2}  \\
\hline \\
\alpha _{s } & \alpha _{s+1} & \alpha _{s+2} & \ldots & \alpha _{k+s-2}  & \alpha _{k+s-1} \\
\end{array}\right). $$

By  $ \eqref{eq 7.11} $ we get
 \begin{equation}
   \label{eq 7.12}
 \sum_{i=0}^{s-2}\Gamma _{i}^{s\times (k-1)}=
 \sum_{i=0}^{s-2} 3\cdot2^{2(i-1)}= 2^{2s-4}
  \end{equation} \\
which implies 
    \begin{equation}
   \label{eq 7.13}
   \Gamma _{s-1}^{s\times (k-1)} +  \Gamma _{s}^{s\times (k-1)}= 2^{k+s-2}-2^{2s-4}.
    \end{equation} \\
   From  Lemma \ref{lem 4.6} it  follows with $ k\longrightarrow k-1 $ and q = 1
\begin{equation}
\label{eq 7.14}
\quad \int_{\mathbb{P}} h_{s,k-1}(t)dt  = \int_{\mathbb{P}} 2^{k+s -1-r( D_{s\times (k-1)}(t))}dt =
  2\cdot\sum_{i=0}^s \Gamma _{i}^{s\times (k-1)}\cdot2^{-i}.
\end{equation}
 On the other hand
 \begin{equation}
 \label{eq 7.15}
  \int_{\mathbb{P}} h_{s,k-1}(t) dt  = 
   Card\left\{(Y,Z),degY\leq k-2,degZ\leq s-1\mid Y\cdot Z = 0\right\}
     =  2^{k-1} +2^s -1. 
\end{equation}
From  \eqref{eqn 7.14} and \eqref{eq 7.15} we obtain
\begin{equation}
\label{eq 7.16}
\sum_{i=0}^s \Gamma _{i}^{s\times(k-1)}\cdot2^{-i}= 2^{k-2}+ 2^{s-1}-2^{-1}.
\end{equation}
In view of \eqref{eq 7.11}
\begin{equation}
\label{eq 7.17}
\sum_{i=0}^{s-2} \Gamma _{i}^{s\times(k-1)}\cdot2^{-i}=
\sum_{i=0}^{s-2}3\cdot2^{2(i-1)}\cdot2^{-i}= 3\cdot2^{s-3} - 2^{-1}.
\end{equation}
From  \eqref{eq 7.16} and \eqref{eq 7.17}
\begin{equation}
\label{eq 7.18}
2\Gamma_{s-1}^{s\times (k-1)} + \Gamma_{s}^{s\times (k-1)}= 2^{k+s-2} +2^{2s-3}.
 \end{equation}
 By    \eqref{eq 7.11}, \eqref{eq 7.13}\; and\;  \eqref{eq 7.18} we deduce  for $ k\geq s+1 $ \\
  \[\Gamma_{i}^{s\times (k-1)} = \left\{\begin{array}{ccc}
             1 & if & i= 0, \\
             3\cdot2^{2(i-1)}  &   if & 1\leq i\leq s-1, \\
   2^{(k+s-2)} - 2^{2s-2}  & if & i = s.
             \end{array}\right.\] \qed


\section{Proof of  Theorem 3.3}
\label{sec 8}
The proof of Theorem 3.3 is based on the relation between
Theorem \ref{thm 3.1} and the following Lemmas.\\
\begin{lem}
\label{lem 8.1}
Let $ t \in \mathbb{P}, \; ( s\leq  k) $ then  \\
\small
\begin{equation}
\label{eq 8.1}
\Gamma _{i}^{s\times k}= \begin{cases}
 {}^{\#}\Big (\begin{array}{c| c}
       0  & 0 \\
       \hline
       0 & 0
\end{array}\Big)_{ \mathbb{P}/\mathbb{P}_{k+s-1}}
& \text{if }  i = 0,  \\
 {}^{\#}\Big (\begin{array}{c| c}
        1 & 1 \\
       \hline
       1 & 1
\end{array}\Big)_{ \mathbb{P}/\mathbb{P}_{k+s-1}}
+{}^{\#}\Big (\begin{array}{c| c}
        0 & 0 \\
       \hline
       0 & 0
\end{array}\Big)_{ \mathbb{P}/\mathbb{P}_{k+s-1}}
& \text{if  }   i = 1, \\
 {}^{\#}\Big (\begin{array}{c| c}
        i & i \\
       \hline
       i & i
\end{array}\Big)_{ \mathbb{P}/\mathbb{P}_{k+s-1}}
+{}^{\#}\Big (\begin{array}{c| c}
        i-1 & i-1 \\
       \hline
       i-1 & i-1
\end{array}\Big)_{ \mathbb{P}/\mathbb{P}_{k+s-1}}
+{}^{\#}\Big (\begin{array}{c| c}
        i-2 & i-1  \\
       \hline
       i-1 & i
\end{array}\Big)_{ \mathbb{P}/\mathbb{P}_{k+s-1}}
 & \text{if  } 2\leq i\leq s-1, \\
  {}^{\#}\Big (\begin{array}{c| c}
        s - 1 & s - 1 \\
       \hline
       s -1 & s -1
\end{array}\Big)_{ \mathbb{P}/\mathbb{P}_{k+s-1}}
+{}^{\#}\Big (\begin{array}{c| c}
        s-1 & s-1 \\
       \hline
       s & s
\end{array}\Big)_{ \mathbb{P}/\mathbb{P}_{k+s-1}}
+{}^{\#}\Big (\begin{array}{c| c}
        s-2 & s-1  \\
       \hline
       s-1 & s
\end{array}\Big)_{ \mathbb{P}/\mathbb{P}_{k+s-1}}
 & \text{if } i = s
\end{cases}
\end{equation}
and \\
\begin{equation}
\label{eq 8.2}
2\cdot\Gamma _{i-1}^{s\times (k-1)}= \begin{cases}
2\cdot{}^{\#}\Big (\begin{array}{c| c}
       0  & 0 \\
       \hline
       0 & 0
\end{array}\Big)_{ \mathbb{P}/\mathbb{P}_{k+s-1}}
& \text{if }  i = 1,  \\
2\cdot{}^{\#}\Big (\begin{array}{c| c}
      i-1 &i -1 \\
       \hline
        i-1 & i-1
\end{array}\Big)_{ \mathbb{P}/\mathbb{P}_{k+s-1}}
+{}^{\#}\Big (\begin{array}{c| c}
        i-2   & i-1\\
       \hline
       i-1 & i
\end{array}\Big)_{ \mathbb{P}/\mathbb{P}_{k+s-1}}
& \text{if  }   2\leq i\leq s-1, \\
2\cdot {}^{\#}\Big (\begin{array}{c| c}
        s -1 & s-1\\
       \hline
       s-1 & s-1
\end{array}\Big)_{ \mathbb{P}/\mathbb{P}_{k+s-1}}
+{}^{\#}\Big (\begin{array}{c| c}
        s-2 & s -1 \\
       \hline
       s -1 & s
\end{array}\Big)_{ \mathbb{P}/\mathbb{P}_{k+s-1}}
& \text{if }    i = s,  \\
{}^{\#}\Big (\begin{array}{c| c}
        s -1 & s-1  \\
       \hline
       s & s
\end{array}\Big)_{ \mathbb{P}/\mathbb{P}_{k+s-1}}
 & \text{if  }  i = s + 1.
\end{cases}
\end{equation}
\end{lem}
\begin{proof}
It suffices to reproduce carefully the proof of Lemma  \ref{lem 5.5}.
\end{proof}

\begin{lem}
\label{lem 8.2}
Let $ t \in \mathbb{P}, \;  ( s\leq  k) $ then \\

\begin{equation}
\label{eq 8.3}
\Gamma _{i}^{s\times k} - 2\cdot\Gamma_{i-1}^{s\times(k-1)}  =
 \begin{cases}
 {}^{\#}\Big (\begin{array}{c| c}
        1 & 1 \\
       \hline
       1 & 1
\end{array}\Big)_{ \mathbb{P}/\mathbb{P}_{k+s-1}}
- {}^{\#}\Big (\begin{array}{c| c}
        0 & 0 \\
       \hline
       0 & 0
\end{array}\Big)_{ \mathbb{P}/\mathbb{P}_{k+s-1}}
& \text{if  }   i = 1, \\
 {}^{\#}\Big (\begin{array}{c| c}
        i & i \\
       \hline
       i & i
\end{array}\Big)_{ \mathbb{P}/\mathbb{P}_{k+s-1}}
-{}^{\#}\Big (\begin{array}{c| c}
        i-1 & i-1 \\
       \hline
       i-1 & i-1
\end{array}\Big)_{ \mathbb{P}/\mathbb{P}_{k+s-1}}
& \text{if     }   2\leq i\leq s-1,  \\
{}^{\#}\Big (\begin{array}{c| c}
        s -1 & s -1  \\
       \hline
       s  & s
\end{array}\Big)_{ \mathbb{P}/\mathbb{P}_{k+s-1}}
- {}^{\#}\Big (\begin{array}{c| c}
        s - 1 & s - 1 \\
       \hline
       s -1 & s -1
\end{array}\Big)_{ \mathbb{P}/\mathbb{P}_{k+s-1}}
& \text{if    }  i = s .
\end{cases}
\end{equation}
\end{lem}
\begin{proof}
Follows from \eqref{eq 8.1} and \eqref{eq 8.2}.
\end{proof}
\begin{lem}
\label{lem 8.3}We have
\begin{equation}
\label{eq 8.4}
\Gamma _{i}^{s\times k} - 2\cdot\Gamma_{i-1}^{s\times(k-1)}  =
 \begin{cases}
 1 & \text{if  } i = 1,\\
 3\cdot2^{2i-3}  &  \text{if     }         2\leq i\leq s-1,\\
 2^{k+s-1}- 5\cdot2^{2s-3}  &    \text{if       }  i = s.
\end{cases}
\end{equation}
\end{lem}
\begin{proof}
The assertion follows immediately from Theorem \ref{thm 3.1}.
\end{proof}
Theorem \ref{thm 3.1} and Lemmas  \ref{lem 8.1}, \ref{lem 8.2}, \ref{lem 8.3}  give the following relations

\footnotesize

\begin{align}
 {}^{\#}\Big (\begin{array}{c| c}
        1 & 1 \\
       \hline
       1 & 1
\end{array}\Big)_{ \mathbb{P}/\mathbb{P}_{k+s-1}}
- {}^{\#}\Big (\begin{array}{c| c}
        0 & 0 \\
       \hline
       0 & 0
\end{array}\Big)_{ \mathbb{P}/\mathbb{P}_{k+s-1}}  & = 1, \label{eq 8.5}\\
 {}^{\#}\Big (\begin{array}{c| c}
        i & i \\
       \hline
       i & i
\end{array}\Big)_{ \mathbb{P}/\mathbb{P}_{k+s-1}}
-{}^{\#}\Big (\begin{array}{c| c}
        i-1 & i-1 \\
       \hline
       i-1 & i-1
\end{array}\Big)_{ \mathbb{P}/\mathbb{P}_{k+s-1}} & = 3\cdot2^{2i-3}\quad \text{if}\quad  2\leq i\leq s-1\label{eq 8.6} \\
{}^{\#}\Big (\begin{array}{c| c}
        s -1 & s -1  \\
       \hline
       s  & s
\end{array}\Big)_{ \mathbb{P}/\mathbb{P}_{k+s-1}}
- {}^{\#}\Big (\begin{array}{c| c}
        s - 1 & s - 1 \\
       \hline
       s -1 & s -1
\end{array}\Big)_{ \mathbb{P}/\mathbb{P}_{k+s-1}} & = 2^{k+s-1} - 5\cdot2^{2s-3}, \label{eq 8.7}\\
 {}^{\#}\Big (\begin{array}{c| c}
        i & i \\
       \hline
       i & i
\end{array}\Big)_{ \mathbb{P}/\mathbb{P}_{k+s-1}}
+{}^{\#}\Big (\begin{array}{c| c}
        i-1 & i-1 \\
       \hline
       i-1 & i-1
\end{array}\Big)_{ \mathbb{P}/\mathbb{P}_{k+s-1}}
+{}^{\#}\Big (\begin{array}{c| c}
        i-2 & i-1  \\
       \hline
       i-1 & i
\end{array}\Big)_{ \mathbb{P}/\mathbb{P}_{k+s-1}} & = 3\cdot2^{2i-2}\quad \text{if} \quad 1\leq i\leq s-1, \label{eq 8.8}\\
{}^{\#}\Big (\begin{array}{c| c}
        s - 1 & s - 1 \\
       \hline
       s -1 & s -1
\end{array}\Big)_{ \mathbb{P}/\mathbb{P}_{k+s-1}}
+{}^{\#}\Big (\begin{array}{c| c}
        s-1 & s-1 \\
       \hline
       s & s
\end{array}\Big)_{ \mathbb{P}/\mathbb{P}_{k+s-1}}
+{}^{\#}\Big (\begin{array}{c| c}
        s-2 & s-1  \\
       \hline
       s-1 & s
\end{array}\Big)_{ \mathbb{P}/\mathbb{P}_{k+s-1}} & = 2^{k+s-1} - 2^{2s-2}. \label{eq 8.9}
\end{align}\vspace{0.3 cm}\\
\normalsize
Now we proceed as follows: \vspace{0.3 cm}\\
 First we deduce by \eqref{eq 8.5} and \eqref{eq 8.6}\vspace{0.3 cm}\\
 \begin{equation} 
 \label{eq 8.10}
   {}^{\#}\Big (\begin{array}{c| c}
        i & i \\
       \hline
       i & i
\end{array}\Big)_{ \mathbb{P}/\mathbb{P}_{k+s-1}}  = 2^{2i-1} 
  \quad   if \; 1\leq i\leq s - 1.  
  \end{equation}\vspace{0.3 cm}\\
Secondly by \eqref{eq 8.8}  and \eqref{eq 8.10}  \vspace{0.3 cm}\\
 \begin{equation*}   
   {}^{\#}\Big (\begin{array}{c| c}
     i-2 & i-1  \\
       \hline
       i-1 & i
\end{array}\Big)_{ \mathbb{P}/\mathbb{P}_{k+s-1}}  =  2^{2i - 3}  \quad   if \; 2\leq i\leq s -1.
  \end{equation*}\vspace{0.3 cm}\\
 Thirdly by  \eqref{eq 8.7}  and  \eqref{eq 8.10} with i = s - 1  \vspace{0.3 cm}
   \begin{equation}
   \label{eq 8.11}    
             {}^{\#}\Big (\begin{array}{c| c}
        s-1 & s-1 \\
       \hline
       s & s
\end{array}\Big)_{ \mathbb{P}/\mathbb{P}_{k+s-1}}  = 2^{k+s-1} - 2^{2s-1}. 
  \end{equation}\vspace{0.3 cm}\\
Last by  \eqref{eq 8.9},  \eqref{eq 8.11} and  \eqref{eq 8.10} with i = s - 1 \vspace{0.3 cm}
 \begin{equation*}  
{}^{\#}\Big (\begin{array}{c| c}
        s-2 & s-1  \\
       \hline
       s-1 & s
\end{array}\Big)_{ \mathbb{P}/\mathbb{P}_{k+s-1}}  =  2^{2s - 3}. 
  \end{equation*}
\section{Proofs of  Theorems 3.4,3.5,3.6,3.7 and 3.8}
\label{sec 9}
\subsection{Proof of Theorem 3.4}
\label{subsec 9.1} Apply   Lemma \ref{lem 4.6} and note that 
\begin{align*}
\int_{\mathbb{P}}h^q(t)dt  = 
\int_{\mathbb{P}}\big[\sum_{deg Y \leq k-1}\sum_{deg Z\leq s-1}E(tYZ)\big]^q dt & =  \\
\int_{\mathbb{P}}\sum_{Y_{1}}\sum_{Z_{1}}\sum_{Y_{2}}\sum_{Z_{2}}\cdots
\sum_{Y_{q}}\sum_{Z_{q}}E(tY_{1}Z_{1})\cdots E(tY_{q}Z_{q})dt  & =  \\
\sum_{Y_{1}}\sum_{Z_{1}}\sum_{Y_{2}}\sum_{Z_{2}}\cdots
\sum_{Y_{q}}\sum_{Z_{q}}\int_{\mathbb{P}}E(t(Y_{1}Z_{1} + \cdots + Y_{q}Z_{q})) dt & =  R
\end{align*}
\qed
\subsection{Proof of Theorem 3.5}
\label{subsec 9.2} Follows from  Lemma \ref{lem 4.10} and  Lemma \ref{lem 5.6}. \qed
\subsection{Proof of Theorem 3.6}
\label{subsec 9.3} Follows from  Lemma \ref{lem 6.2} and  Lemma \ref{lem 6.4}. \qed
\subsection{Proof of Theorem 3.7}
\label{subsec 9.4} Immediately  from  Lemma \ref{lem 6.9} . \qed
\subsection{Proof of Theorem 3.8}
\label{subsec 9.5}
Let  $ 2\leq m\leq k-2 .$ Then, by  Lemma \ref{lem 6.8} and Theorem \ref{thm 3.1} with s = m+1\\

\begin{align*}
  \Gamma _{0}^{\Big[\substack{1 \\ 1+m }\Big] \times k}       & = 1, \\
  \Gamma _{1}^{\Big[\substack{1 \\ 1+m }\Big] \times k}   & = \sigma _{0,1} + \sigma _{1,1} = (2^k - 1) + 2\cdot3 = 2^k +5, \\
  \Gamma _{i}^{\Big[\substack{1 \\ 1+m }\Big] \times k}      & = \sigma _{i-1,i} + \sigma _{i,i} =(2^k-2^{i-1})\cdot \Gamma _{i-1}^{(1+m)\times k}
 + 2^{i}\cdot\Gamma _{i}^{(1+m)\times k}  \\
 & = (2^k-2^{i-1})\cdot3\cdot2^{2(i-2)}  + 2^{i}\cdot3\cdot2^{2(i-1)}) = 3\cdot2^{k+ 2i -4 } + 21\cdot2^{3i -5} && \text{if   } 2\leq i\leq  m, \\
   \Gamma _{m+1}^{\Big[\substack{1 \\ 1+m }\Big] \times k}          & = \sigma _{m,m+1} + \sigma _{m+1,m+1} =
  (2^k-2^{m})\cdot \Gamma _{i-1}^{m\times k} +  2^{m+1}\cdot\Gamma _{m+1}^{(1+m)\times k} \\
  & = (2^k-2^{m})\cdot3\cdot2^{2(m-1)}  + 2^{1+m}\cdot(2^{k+m} - 2^{2m}) = 11\cdot[2^{k+2m-2} - 2^{3m -2}],   \\
  \Gamma _{m+2}^{\Big[\substack{1 \\ 1+m }\Big] \times k}   & = \sigma _{m+1,m+2} = (2^k - 2^{1+m})\cdot\Gamma _{m+1}^{(1+m)\times k} =  (2^k - 2^{1+m})\cdot(2^{k+m} - 2^{2m}) \\
   & = 2^{2k+m}  - 3\cdot2^{k+2m} +2^{3m+1}. 
\end{align*}
The proof in the case $ 3\leq  k \leq 1+m $ is similar.\qed

\section{Proofs of  Theorem  3.9 , Corollary 3.10 and Theorem 3.11}
\label{sec 10}
\subsection{ $ \Gamma_{i}^{\left[n\atop 1+m\right]\times k} $ written as a linear  combination of the  $\Gamma _{i-j}^{(1+m)\times k}$ for $ 0\leq j\leq n. $}
\label{subsec 10.1}
  \begin{lem}
\label{lem 10.1}
[See Definitions \ref{defn 1.1},   \ref{defn 1.2}].              
The number $ \Gamma _{i}^{\Big[\substack{n \\ 1+m }\Big] \times k}, $ expressed as a linear combination of the $ \Gamma _{i-j}^{(1+m)\times k} $
for $ j = 0,1,\ldots n, $ is given by the following recurrence formula 
\begin{equation}
\label{eq 10.1}
 \Gamma _{i}^{\Big[\substack{n \\ 1+m }\Big] \times k}= 
 \sum_{j= 0}^{n}\Big[2^{(n-j)\cdot(i-j)} a_{j}^{(n)}\prod_{l=1}^{j}(2^{k}- 2^{i-l})\Big]\cdot
 \Gamma _{i-j}^{(1+m)\times k}, n= 1,2,\ldots \quad for \quad 0\leq i\leq inf(k,n+m +1),
\end{equation} 
 where $ a_{j}^{(n)} $ satisfies the linear recurrence relation 
 \begin{equation}
 \label{eq 10.2}
 a_{j}^{(n)} = 2^{j}\cdot a_{j}^{(n-1)} + a_{j-1}^{(n-1)},\quad n = 2,3,4,\ldots       \quad for\quad 1\leq j\leq n-1. 
\end{equation}
 We set  
   \begin{align*}
   a_{0}^{(n)} & =   a_{n}^{(n)} = 1  \\
 and \quad  \Gamma _{i-j}^{(1+m)\times k} & = 0 \quad if \quad  i-j \notin \{0,1,2,\ldots, inf(k,1+m)\}.
   \end{align*}
   We can then compute successively
      $\quad a_{j}^{(n)},\quad  n\geq j \quad for \quad  j= 1,2,\ldots \; and\; we\; obtain $
  \begin{align*}  
 a_{1}^{(n)} & = 2^n -1\quad for\quad n\geq 1, \\
 a_{2}^{(n)} & = {2^{2n-1}-3\cdot2^{n-1} + 1\over 3} \quad  for \quad n\geq 2, \\
 a_{j}^{(n)} & =  2^{nj-j^2} + \sum_{l=0}^{n-j-1} a_{j-1}^{(n-1-l)}\cdot(2^j)^l \quad  for \quad n\geq j.  
    \end{align*}  
   \end{lem}
\begin{proof}
We prove Lemma \ref{lem 10.1} by induction on n. \\
A similar proof  of  Lemma \ref{lem 6.8} gives the following generalization
\begin{equation}
\label{eq 10.3}
  \Gamma _{i}^{\Big[\substack{n \\ 1+m }\Big] \times k}= 2^{i}\cdot \Gamma _{i}^{\Big[\substack{n-1 \\ 1+m }\Big] \times k}    
    +      (2^{k}-2^{i-1})\cdot \Gamma _{i-1}^{\Big[\substack{n-1 \\ 1+m }\Big] \times k}\quad for \quad 0\leq i\leq inf(k,n+m +1).
    \end{equation}

 From Lemma \ref{lem 6.8} and \eqref{eq 10.3} with n=2 we get respectively for n =1 and n = 2 \\
\begin{align}
  \Gamma _{i}^{\Big[\substack{1 \\ 1+m }\Big] \times k} &  =  
  2^{i}\Gamma _{i}^{(1+m)\times k} + ( 2^{k} -  2^{i-1})\cdot\Gamma _{i-1}^{(1+m)\times k} 
   \quad  for \quad  0\leq i\leq inf(k,2+m) \label{eq 10.4} \\
  \Gamma _{i}^{\Big[\substack{ 2 \\ 1+m }\Big] \times k} 
& =   2^{i}\Gamma _{i}^{\left[\stackrel{1}{1+m}\right]\times k} +
 ( 2^{k} -  2^{i-1})\cdot \Gamma _{i-1}^{\left[\stackrel{1}{1+m}\right]\times k}\quad for \quad  0\leq i\leq inf(k,3+m). \label{eq 10.5}
 \end{align}
From (10.4) and (10.5) we deduce \\
\begin{align}
\Gamma _{i}^{\Big[\substack{2 \\ 1+m }\Big] \times k}
& = 2^{2i}\Gamma _{i}^{(1+m)\times k}+
     3\cdot2^{i-1}(2^{k}-2^{i-1})\cdot\Gamma _{i-1}^{(1+m)\times k} \label{eq 10.6}\\
       &  +(2^{k}-2^{i-1})(2^{k}-2^{i-2})\cdot
      \Gamma _{i-2}^{(1+m)\times k}  \quad for\quad 0\leq i\leq inf(k,3+m). \nonumber
\end{align}
Hence by (10.4) and (10.6) the formula (10.1) holds for n = 1 and n = 2. \\
Assume now that (10.1) holds for the number n-1, that is
 \begin{align}    
   2^{i}\cdot \Gamma _{i}^{\Big[\substack{n-1 \\ 1+m }\Big] \times k} =    
 \sum_{j= 0}^{n-1}2^{(n-1-j)\cdot(i-j)}2^{i} a_{j}^{(n-1)}\prod_{l=1}^{j}(2^{k}- 2^{i-l})\cdot
 \Gamma _{i-j}^{(1+m)\times k},\label{eq 10.7}
 \end{align}
\begin{align}
 (2^{k}-2^{i-1})\cdot \Gamma _{i-1}^{\Big[\substack{n-1 \\ 1+m }\Big] \times k} & =
  (2^{k}-2^{i-1})\big[ \sum_{j= 0}^{n-1}2^{(n-1-j)\cdot(i-1-j)} a_{j}^{(n-1)}\prod_{l=1}^{j}(2^{k}- 2^{i-1-l})\cdot
 \Gamma _{i-1-j}^{(1+m)\times k}\big] \label{eq 10.8} \\ 
 & =
  \sum_{j= 0}^{n-1}2^{(n-1-j)\cdot(i-1-j)} a_{j}^{(n-1)}\prod_{l=1}^{j+1}(2^{k}- 2^{i-l})\cdot
 \Gamma _{i-(j+1)}^{(1+m)\times k} \nonumber \\
 & = \sum_{j= 1}^{n}2^{(n-j)\cdot(i-j)} a_{j-1}^{(n-1)}\prod_{l=1}^{j}(2^{k}- 2^{i-l})\cdot \Gamma _{i-j}^{(1+m)\times k}. \nonumber
\end{align}
From \eqref{eq 10.7}, \eqref{eq 10.8}  ,  and \eqref{eq 10.2} it follows
 \begin{align}
 \Gamma _{i}^{\Big[\substack{n \\ 1+m }\Big] \times k}
 & = 2^{ni}\Gamma _{i}^{(1+m)\times k} 
+ \sum_{j= 1}^{n-1}2^{(n-j)\cdot(i-j)}(2^{j}a_{j}^{(n-1)}  +
 a_{j-1}^{(n-1)})\prod_{l=1}^{j}(2^{k}- 2^{i-l})\cdot \Gamma _{i-j}^{(1+m)\times k} \label{eq 10.9}\\
& + a_{n-1}^{(n-1)}\prod_{l=1}^{n}(2^{k}- 2^{i-l})\Gamma _{i-n}^{(1+m)\times k} \nonumber  \\
&  =  \sum_{j= 0}^{n}2^{(n-j)\cdot(i-j)} a_{j}^{(n)}\prod_{l=1}^{j}(2^{k}- 2^{i-l})\cdot
 \Gamma _{i-j}^{(1+m)\times k}. \nonumber
\end{align}

\end{proof}
\subsection{Computation of $ a_{i}^{n} $ for  $1 \leq i \leq n-1 $}
\label{subsec 10.2}
  \begin{lem}
\label{lem 10.2}
We have for $ 1\leq i\leq n-1$
\begin{equation*}
 a_{i}^{(n)} = \sum_{s =0}^{i-1} (-1)^{s}\prod_{l=0}^{i-(s+1)}{2^{n+1}- 2^{l}\over 2^{i-s}-2^{l}}\cdot2^{s(n-i) +{s(s+1)\over 2}}
+ (-1)^{i}\cdot2^{in - {i(i-1)\over 2}}. 
\end{equation*}
\end{lem}
\begin{proof}
(See Definition \ref{defn 1.2})
We recall that $ \Gamma _{i}^{\Big[\substack{n \\ 1 }\Big] \times k}$ denotes  the number 
of $(n+1)\times k $ matrices over $ \mathbb{F}_{2} $ of rank i.\\
By formula \eqref{eq 10.1} with m = 0, we obtain\\
\begin{equation}
\label{eq 10.10}
 \Gamma _{i}^{\Big[\substack{n \\ 1 }\Big] \times k}= 
 \sum_{j= 0}^{n}2^{(n-j)\cdot(i-j)} a_{j}^{(n)}\prod_{l=1}^{j}(2^{k}- 2^{i-l})\cdot
 \Gamma _{i-j}^{ 1\times k}, n= 1,2,\ldots \quad for \quad 0\leq i\leq inf(k,n +1). 
\end{equation}
We have obviously \begin{equation}
\label{eq 10.11}
  \Gamma _{i-j}^{ 1\times k}  = \begin{cases}
 1 & \text{if      }   i-j = 0, \\
  2^{k} -1 & \text{if      }   i-j = 1.
    \end{cases}
\end{equation}
From \eqref{eq 10.10} and  \eqref{eq 10.11} we deduce
\begin{align}
\Gamma _{i}^{\Big[\substack{n \\ 1 }\Big] \times k} &  = \sum_{j= i-1}^{i}2^{(n-j)\cdot(i-j)} a_{j}^{(n)}\prod_{l=1}^{j}(2^{k}- 2^{i-l})\cdot
 \Gamma _{i-j}^{ 1\times k} \label{eq 10.12}\\
 & = a_{i-1}^{(n)}\cdot2^{n-(i-1)}\prod_{l=1}^{i}(2^{k}- 2^{i-l}) +  a_{i}^{(n)}\cdot\prod_{l=1}^{i}(2^{k}- 2^{i-l})\nonumber  \\
 & = \prod_{l=1}^{i}(2^{k}- 2^{i-l})[ a_{i}^{(n)} +  a_{i-1}^{(n)}\cdot2^{n-(i-1)}]= 
  \prod_{l=0}^{i-1}(2^{k}- 2^{l})[ a_{i}^{(n)} +  2^{n-(i-1)}\cdot a_{i-1}^{(n)} ].  \nonumber 
\end{align} 
On the other hand by George Landsberg [3] we have
\begin{equation}
\label{eq 10.13}
\Gamma _{i}^{\Big[\substack{n \\ 1 }\Big] \times k}=
 \prod_{l = 0}^{i-1}{ (2^{n+1} -2^{l})(2^{k}-2^{l}) \over (2^{i}- 2^{l})}. 
\end{equation} 
 Hence by \eqref{eq 10.12} and \eqref{eq 10.13} we have the formula
 \begin{equation}
 \label{eq 10.14}
 a_{i}^{(n)} +  2^{n-(i-1)}\cdot a_{i-1}^{(n)} =  \prod_{l = 0}^{i-1}{ 2^{n+1} -2^{l} \over 2^{i}- 2^{l}}. 
 \end{equation}
From  \eqref{eq 10.14} we deduce 
\begin{equation}
\label{eq 10.15}
           \sum_{s = 0}^{i-1} (-1)^{s}\cdot2^{s(n-i) + {s(s+1)\over 2}}\left[ a_{i-s}^{(n)} +  2^{n-(i-(s+1))}\cdot a_{i-(s+1)}^{(n)}\right] = 
         \sum_{s = 0}^{i-1}\prod_{l = 0}^{i-(s+1)}{ 2^{n+1} -2^{l} \over 2^{i-s}- 2^{l}} (-1)^{s}\cdot2^{s(n-i) + {s(s+1)\over 2}}.
  \end{equation}
  Using  \eqref{eq 10.15} we get after some simplifications
  \begin{equation*}
   a_{i}^{(n)}  - (-1)^{i}a_{0}^{(n)}2^{in - {i(i-1)\over 2}}= 
  \sum_{s = 0}^{i-1}\prod_{l = 0}^{i-(s+1)}{ 2^{n+1} -2^{l} \over 2^{i-s}- 2^{l}} (-1)^{s}\cdot2^{s(n-i) + {s(s+1)\over 2}}. 
\end{equation*}
 \end{proof}
\subsection{ Proof of Theorem \ref{thm 3.9}}
\label{subsec 10.3}
  Theorem \ref{thm 3.9} follows from Lemma \ref{lem 10.1} and Lemma  \ref{lem 10.2}.
\subsection{ Proof of Corollary \ref{cor 3.10}}
\label{subsec 10.4}
 The assertions follows from Theorem \ref{thm 3.9} by some simple calculations.
 \subsection{ Proof of Theorem  \ref{thm 3.11}}
\label{subsec 10.5}
The proof of Theorem \ref{thm 3.11} is just a generalization of the proof of Lemma \ref{lem 6.9}.


\begin{thebibliography}{99}
\bibitem{Daykin}  Daykin David E,  {\textit{Distribution of Bordered Persymmetric Matrices in a finite field}}
{J. reine angew. Math}, {\bf 203}, 47-54 (1960). 
\bibitem{Hayes}Hayes , D.R,  {\textit{The expression of a polynomial as a sum of three irreducibles}}
{Acta Arith.} {\bf 11}, 461-488 (1966).
\bibitem{Landsberg}Landsberg, G {\textit{Ueber eine Anzahlbestimmung und eine damit zusammenhangende Reihe}},
 {J. reine angew. Math.}, {\bf 111}, 87-88 (1893). 
\end{thebibliography}
\end{document}